\documentclass[aip,
 amsmath,amssymb,
 reprint,%
]{revtex4-2}

\usepackage{graphicx}% Include figure files
\usepackage{dcolumn}% Align table columns on decimal point
\usepackage{bm}% bold math
\usepackage[utf8]{inputenc}
\usepackage[T1]{fontenc}
\usepackage{mathptmx,xcolor,upgreek,enumitem,amsthm,yfonts}
\usepackage{comment} 
\usepackage{etoolbox}
\usepackage{overpic}
\usepackage[colorlinks=true,linkcolor=red!60!black,citecolor=green!60!black]{hyperref} %para PDFLaTeX
\usepackage{makecell}

\newtheorem{theorem}{Theorem}[section]
\newtheorem{definition}[theorem]{Definition}

\newtheorem{proposition}[theorem]{Proposition}

\newcommand{\R}{\mathbb{R}}

\makeatletter
\def\@email#1#2{%
 \endgroup
 \patchcmd{\titleblock@produce}
  {\frontmatter@RRAPformat}
  {\frontmatter@RRAPformat{\produce@RRAP{*#1\href{mailto:#2}{#2}}}\frontmatter@RRAPformat}
  {}{}
}%
\makeatother
\begin{document}

\preprint{AIP/123-QED}

\title[Rate-induced tracking in a time-dependent environment]{Rate-induced tracking for concave or d-concave transitions in a time-dependent environment with application in ecology}
% Force line breaks with \\
\author{J. Due\~{n}as}
\affiliation{Departamento de Matem\'{a}tica Aplicada, Universidad de Va\-lladolid\\ Paseo Prado de la Magdalena 3-5, 47011 Valladolid, Spain.}
\affiliation{Instituto de Investigaci\'{o}n en Matem\'{a}ticas (IMUVa), Universidad de Valladolid\\ Paseo de Bel\'{e}n S/N, 47011 Valladolid, Spain}

 \author{I.P. Longo}%
\affiliation{Imperial College London, Department of Mathematics, DynamIC, 637 Huxley Building, 180 Queen's Gate South
Kensington Campus, London, SW7 2AZ, United Kingdom.}

\author{R. Obaya}
\affiliation{Departamento de Matem\'{a}tica Aplicada, Universidad de Va\-lladolid\\ Paseo Prado de la Magdalena 3-5, 47011 Valladolid, Spain.}
\affiliation{Instituto de Investigaci\'{o}n en Matem\'{a}ticas (IMUVa), Universidad de Valladolid.}

\email{rafael.obaya@uva.es}
\email{iacopo.longo@imperial.ac.uk}
\email{jesus.duenas@uva.es}

\date{\today}% It is always \today, today,

\begin{abstract}
This paper investigates biological models that represent the transition equation from a system in the past to a system in the future.
It is shown that finite-time Lyapunov exponents calculated along a locally pullback attractive solution are efficient indicators (early-warning signals) of the presence of a tipping point.
Precise time-dependent transitions with concave or d-concave variation in the state variable giving rise to scenarios of rate-induced tracking are shown. They are classified depending on the internal dynamics of the set of bounded solutions. Based on this classification, some representative features of these models are investigated by means of a careful numerical analysis.
\end{abstract}

\maketitle

\begin{quotation}
In those transitions where the parameter shift provides intermediate systems with an inappropriate hyperbolic structure for a long period of time, rate-induced tracking can occur.
In these transitions, high rates are beneficial, whereas a tipping point may occur when the transition rate decreases excessively.
In this paper, we analyze this phenomenon in precise concave or d-concave transitions of biological models.
In these examples, we show that finite-time Lyapunov exponents are effective early-warning signals of critical transitions, since they experience an overgrowth when the transition rate is close to a tipping point at which the dynamics changes drastically.
In this formulation, a mechanism of control used as response to warning signal gives rise to safe operating margins in rate-induced tracking problems with time-dependent rate.
\end{quotation}

\section{Introduction}\label{sec:1introduction}
Critical transitions are large and often irreversible changes in the state of a given system in response to small and slow changes in the input.
Remarkable phenomena in climate,\cite{aajqw,awvc} biology\cite{hill} and ecology\cite{atw,scheffer2} have been described under the headline of tipping.

Since Ashwin et al.,\cite{aspw} it has become common to mathematically formulate the problem as a time-dependent transition connecting a system in the past to a system in the future.\cite{alas,altw,awvc,lno1,lno2,lno3,dno3,kuehnlongo}
It is also assumed that the long-term dynamics of past and future systems are concentrated in an appropriate local or global attractor.
In some cases, the evolution of the attractor of the past, corresponding to a pullback attractor of the transition equation, connects to the attractor of the future. This behavior is called tracking. Hovewer, in other cases the pullback attractor loses stability, changes its shape or even dissapears. Then, the connection between  the past and the future attractors fails.
This behaviour is called tipping. Within the framework of transition equations, scenarios of rate-induced tipping, phase-induced tipping and size-induced tipping have been extensively investigated in the literature.\cite{aspw,alas,altw,awvc,lno1,lno2,lno3,dno3}
In this paper, we consider models given by the differential equation
\begin{equation}\label{eq:x'=f(t,x,gamma)}
    x'=f(t,x,\gamma)\,,
\end{equation}
where $f\colon\mathbb{R}\times\mathbb{R}\times\mathbb{R}\rightarrow\mathbb{R}$ is smooth and satisfies one of the following two properties with respect to the state variable: it is concave or it has a concave partial derivative. (In the second case, we say that $f$ is d-concave for short).\cite{dno1}
The past and the future systems are respectively given by \eqref{eq:x'=f(t,x,gamma)} with $\gamma=\gamma_-$ and $\gamma=\gamma_+$. Moreover, a transition equation between past and future system is given by
\begin{equation}\label{eq:x'=f(t,x,Gamma(t))}
    x'=f(t,x,\Gamma(t))\,,
\end{equation}
where $\Gamma\colon\mathbb{R}\rightarrow\mathbb{R}$ is continuous and $\lim_{t\rightarrow\pm\infty}\Gamma(t)=\gamma_\pm$.
Typically, one is interested in a family of transition equations, 
\begin{equation}\label{eq:x'=f(t,x,Gammac(t))}
    x'=f(t,x,\Gamma^c(t))\,,
\end{equation}
indexed by a parameter $c$, whose variation may sometimes lead to a critical transition, i.e.~a change from a tracking to a tipping dynamical scenario. Usually, it can be rigorously shown that the mathematical mechanism behind such critical transition is a nonautonomous saddle-node bifurcation.

%A critical transition for \eqref{eq:x'=f(t,x,Gamma(t))} may appear as a parameter $c$ varies in a family of transition equations 
%\begin{equation}\label{eq:x'=f(t,x,Gammac(t))} 
%x'=f(t,x,\Gamma^c(t))\,,
%\end{equation}
%and, usually, it can be rigorously shown to be the consequence of a nonautonomous saddle-node bifurcation. 

In this paper, we investigate the occurrence of rate-induced tracking in \eqref{eq:x'=f(t,x,Gammac(t))} with  $\Gamma^c(t)=\Gamma(ct)$, i.e.~where the parameter $c$ is the rate at which the transition takes place.
In contrast to rate-induced tipping, where tracking is guaranteed for small rates and tipping occurs beyond a certain critical rate, in rate-induced tracking the opposite scenario holds: there is always tipping for sufficiently small rates and tracking beyond a certain rate.
This concept covers a range of problems for which high rates are beneficial: they prevent the occurrence of tipping. As it is exposed in Section~\ref{sec:3biologicalmodels}, concave and d-concave functions often appear in mathematical models of biology and ecology. All these comments justify the unified treatment of concave and d-concave transition theory posed in Section~\ref{sec:4tippingtracking} of this work.

Early-warning signals (EWS) of critical transitions are of huge importance in many areas of applied sciences, since they allow to make all the possible efforts to avoid a catastrophe.
Recovery rates, which correspond to some extent to the mathematical notion of finite-time Lyapunov exponents (FTLE), have been repeatedly highlighted in the literature on applied sciences as EWS of critical transitions.\cite{lenton1,das,moore,schefferbas,dakos}
Most of these works are devoted to the study of EWS of critical transitions with slow parametric variation.
That is, dynamical features of the intermediate equations \eqref{eq:x'=f(t,x,gamma)} which point out the proximity of a bifurcation point in the parameter $\gamma$ are sought.
Seen from the perspective of the transition equation, that approach implicitly assumes that the variation of $\Gamma$ is so slow that the system stays sufficiently close to each of the elements of the curve of stable solutions of the intermediate equations over a time span which is long enough to reproduce those dynamical features before reaching a bifurcation point.
In particular, the decrease to zero of the recovery rates has been frequently used as an EWS and usually called critical slowing down.

In nonautonomous dynamics, the possible existence of strange nonchaotic attractors \cite{milon,vinograd,johnson1,jager3}
at the bifurcation point of the $\gamma$-parametric family
(where both strictly positive and negative Lyapunov exponents exist)
may prevent the occurrence of such a critical slowing down.
Nonetheless, it has been found that, even in the case of existence of a strange nonchaotic attractor, there exists a set of positive measure in the skewproduct base for which a change of sign of the FTLE takes place, as long as the FTLE are calculated for a sufficiently small integration time.\cite{jager1}

In Section~\ref{sec:5earlywarningsignals} of this paper,
we describe EWS of transitions which do not necessarily have slow parametric variation.
In particular, FTLE are shown to be effective EWS for $c$-parametric families of transition equations \eqref{eq:x'=f(t,x,Gammac(t))}  with $\Gamma^c(t)=\Gamma(ct)$ which exhibit rate-induced tracking, that is, for which high rates $c>0$ (fast transition) are beneficial.
The state of the system during the transition is described by a locally pullback attractive solution. The FTLE along this special trajectory exhibit a change of sign near the critical value of $c$ where a drastic change of asymptotic dynamics appears.
As a matter of fact, a substantial increase of the FTLE (even before reaching positive values) may be already considered as an EWS.

We remark that the skewproduct formulation in Subsection~\ref{subsec:4jointskewproduct} univoquely identifies the trajectory on which this phenomenon takes place.
%In Subsection \ref{subsec:5Crateinduced}, we provide an example on which this new EWS and the critical slowing down previously described in the literature appear for different values of $c$ in the same $c$-parametric family: the critical slowing down appears for $c$ near 0 while the new EWS appears for $c$ around the bifurcation point of the $c$-parametric family of transition equations.
Although our analysis focuses on rate-induced tracking, we stress that the proposed EWS can also be effective for other critical transitions appearing in equations of the type of \eqref{eq:x'=f(t,x,Gammac(t))}.

In Section~\ref{sec:Vrateinducedtrackingphenomena}, we investigate several transition equations for concave and d-concave time-dependent systems where rate-induced tracking occurs.
In particular, Subsection \ref{subsec:5Crateinduced}, deals with a population modelled by the equation
\begin{equation}
    x'=R(t,x)+\gamma\,\phi (t)\,,
\end{equation}
where $R$ is d-concave and 
the migration term $\gamma\,\phi (t)$ is multiple of a strictly positive function $\phi(t)$.
We will work under assumptions that guarantee the existence of a compact interval $I$ such that, if $\gamma\in I$, then there exists a solution representing a good-sized population.
On the other hand, if $\gamma < \inf I$, then the population decreases to extinction.
Some migration phenomena justify the substitution of $\gamma$ by a rate-dependent parameter shift $\Gamma(ct)$ which departs from a value $\gamma_+ \in I$ in the past, decreases monotonically to $\Gamma(0)=\gamma_*<\inf I$, and returns to $\gamma_+$ in the future.
It is known that there exists a value $c_0$ such that, if the rate $c>c_0$, then the population of the transition equation keeps a good size, and that, if $c<c_0$, then the population dramatically decreases to extinction.\cite{dno3} 
A detailed numerical analysis of the main characteristics of these transitions is presented.
We take $\kappa \in (0,1)$ and issue an EWS each time that $\lambda_u (c,T,t)>\kappa\,L$, where $\lambda_u(c,T,t)$ is the FTLE of the upper locally pullback attractive solution $u_c$ on the interval $[t-T,t]$ and $L$ is the Lyapunov exponent of the upper hyperbolic attractive solution of the past.
Under this formulation, we analyze the range of rates $c>c_0$ for which the previous EWS points out the proximity of a tipping point.
Next, for some pairs $(c,\kappa)$, we determine the values of time $t$ and lengths $T$ for which the criteria $\lambda_u (c,T,t)>\kappa\, L$ is satisfied and thus an EWS is issued.
We continue the study in the case where $\Gamma(\Delta(t)\,t)$ has time-dependent transition rate $\Delta(t)$, which is considered as a new parameter shift, assuming that there exist mechanisms able to change its value over time.

In Subsection~\ref{subsec:5Cconcaverateinducedtracking}, we focus our attention on the biological logistic model with time-dependent rate and a migration term of the form 
\begin{equation}\label{eq:variablerateintroduction}
    x'=-\big(x-\Gamma(\Delta(t)\,t)\big)^2+I(t)\,.
\end{equation}
The underlying assumptions will be that that the quadratic equation $x'=-x^2+I(t)$ has an attractor-repeller pair and that $\Gamma(t)$ is asymptotically constant with $\lim_{t \rightarrow \pm \infty}\Gamma(t)=\gamma_\pm$.
Moreover, we will assume that the rate function $\Delta$ evolves from a value $c_1$ in the past to a value $c_2$ in the future with $c_1<\Delta(0)=c_0<c_2$ and such that the transition equation 
\begin{equation}
    x'=-\big(x-\Gamma(ct)^2\big)+I(t)
 \end{equation}
connects the attractor-repeller pair of the past with that of the future when $c\in\{c_1,c_2\}$ but exhibits tipping for $c=c_0$.
Under these conditions, the response of \eqref{eq:variablerateintroduction} depends on the rate of $\Delta$ and a new rate-induced tracking phenomenon can occur.
A similar analysis is carried out in the logistic model 
\begin{equation}
x'=-\big(x-\Gamma(t-\Delta(t))\big)^2+I(t)\,,
\end{equation}
where now $\Delta$ connects two different values of the phase, obtaining again the same previous conclusions.
Finally, in Subsection~\ref{subsec:6Cnoreturnrateinducetracking}, we return to the problems of Subsection \ref{subsec:5Crateinduced} to introduce a control mechanism in the model in order to speed up the transition once an EWS has been issued, showing numerically that a rapid action during a short period of time can avoid the occurrence of tipping.

At the end of the paper, Table~\ref{tab:table3} can be found, which contains a summary of the most commonly used symbols.
\section{Concave and d-concave biological models}\label{sec:3biologicalmodels}
Many features of life on Earth fluctuate over time. A range of phenomena from Earth's rotation to climate variability or the alternation of seasons influences the environments of populations.\cite{renshaw1}
Those populations whose law of evolution (the set of factors affecting their development) is time-dependent can be suitably modeled by nonautonomous equations
\begin{equation}\label{eq:2Ax'=h}
    x'=h(t,x)\,,
\end{equation}
where $t$ represents time and the state variable $x$ represents the size of the population.
From now on, as long as no further conditions of regularity are required, the population growth rate $h\colon\mathbb{R}\times\mathbb{R}\rightarrow\mathbb{R}$ is a continuous function such that its derivative $h_x$ with respect the second variable exists, is continuous, and the restriction of $h$ and $h_{x}$ to $\mathbb{R}\times K$ are bounded and uniformly continuous for any compact $K\subset\mathbb{R}$.
These conditions ensure that the time variation of $h$ is compactifiable,\cite{selltopdyn,shenyi} which is one of the main assumptions of the paper.
This kind of functions include as simplest examples autonomous, periodic, almost periodic and asymptotically almost periodic functions.
%Given the initial datum $(s,x_0)\in\mathbb{R}^2$, we denote by $x(\cdot,s,x_0)\colon(\alpha_{s,x_0},\beta_{s,x_0})\rightarrow\mathbb{R}$, the maximal solution of \eqref{eq:2Ax'=h} with $x(s,s,x_0)=x_0$.
Under this nonautonomous approach, all the model parameters (growth rates, carrying capacities, competitive effects, predation features...) may be time-dependent functions.
%As already mentioned at the beginning of Section~\ref{sec:2preliminaries}, we require the time-variation of the equations to be compactificable in all the models.

Many continuous models of single-species populations in mathematical biology are concave or d-concave,\cite{dno1} that is, the population growth rate $h$ as a function of the population size $x$ either is concave or has concave derivative. For this reason, in this paper, we supply a joint approach to these two types of population models.
In what follows, $h$ is said to be \emph{concave} if
\begin{equation*}
    h(t,\,\alpha x+(1-\alpha)\,y)\geq\alpha \,h(t,x)+(1-\alpha)\,h(t,y)
\end{equation*}
for all $t,x,y\in\mathbb{R}$ and $\alpha\in[0,1]$. If $h_x$ is concave, then $h$ is said to be \emph{d-concave}.\cite{dno1,dno2}
In the case that $h$ is sufficiently regular, it is known that $h$ is concave (resp. d-concave) if and only if $h_{xx}$ (resp. $h_{xxx}$) is a nonpositive function.

To motivate the rest of the paper, we present some concave and d-concave nonautonomous counterparts of classical single-species population equations.\cite{renshaw1,dno3}
%The evolution of the state variable $x$, which represents the size of the population, is governed by the scalar nonautonomous ordinary differential equation $x'=h(t,x)$. As already said, in all the considered cases either $x\mapsto h(t,x)$ or $x\mapsto h_x(t,x)$ is concave for all $t\in\mathbb{R}$. Moreover, the presented models are strictly concave or strictly d-concave in the sense of Subsection~\ref{subsec:2Dstrictconcavityanddconcavity}.
We start from the autonomous logistic equation 
\begin{equation}\label{eq:3autonomousconcavequadratic}
    x'=r\,x\,(1-x/K)\,,
\end{equation}
where $r,K>0$ stand for the intrinsic growth rate of the population and the carrying capacity of the environment, respectively. The simplest example of nonautonomous concave equation is the nonautonomous logistic equation\cite{renshaw1}
\begin{equation}\label{eq:3concavequadratic}
x'=r(t)\,x\,(1-x/K(t))\,.
\end{equation}
It is remarkable that, in this model, the positively bounded-from-below map $t\mapsto K(t)$ does no longer represent the healthy steady population since it is not even a solution of the equation, though the positively bounded-from-below map $t\mapsto r(t)$ retains the meaning that $r$ has in \eqref{eq:3autonomousconcavequadratic}.

Among other available concave models, we find the nonautonomous Gompertz model,\cite{montroll1} which precludes the exponential growth of the population multiplying the right-hand side function by a logarithm instead of a linear term,
\begin{equation}
x'=-r(t)\,x\,\log(x/K(t))\,.
\end{equation}
The nonautonomous continuous Beverton-Holt model,\cite{dossantos1} which has been shown to fit more accurately data of exploited fish populations, is not concave in all the real line but in a neighborhood of the positive halfline, where the biologically significant dynamics takes place: it is given by
\begin{equation}
x'=x\, \left(\frac{1+r(t)}{1+\alpha(t)\,x}-1\right)\,,
\end{equation}
where the positively bounded-from-below map $t\mapsto \alpha(t)$ represents the per capita intraspecific competitive effect as a function of time. The role played by $\alpha$ in this model is related to that of $K$ in the previous ones, since it precludes an infinite growth of the population.

The mathematical modeling of the Allee effect,\cite{courchamp1} which consists on a positive correlation between the size of a population and the population growth rate per individual, provides differential equations given by functions with concave derivative.\cite{dno3} There are several biological mechanisms (easier mate finding, cooperative breeding, cooperative anti-predator behavior, increased foraging efficiency...) responsible for the Allee effect.
The Allee effect can be added to the logistic equation \eqref{eq:3concavequadratic} either in multiplicative form or in additive form. Depending on its features, one of the following multiplicative models\cite{courchamp1,altw} may be preferred to the other one:
\begin{equation}
x'=r(t)\,x\,\left(1-\frac{x}{K(t)}\right)\frac{x-S(t)}{K(t)}\,,
\end{equation}
\begin{equation}\label{eq:3dconcavenoncubicallee}
 x'=r(t)\,x\,\left(1-\frac{x}{K(t)}\right)\frac{x-\mu(t)}{\nu(t)+x}\,,
\end{equation}
where $S$, $\mu$ and $\nu$ are maps which determine the strength of the Allee effect, $K(t)+S(t)\geq 0$ for all $t\in\mathbb{R}$, and $\nu$ and $\nu+\mu$ are positively bounded-from-below.

If the triggering mechanism of the Allee effect is related to predation, an additive form, consisting on adding a Holling type II functional response term to the logistic equation, would be preferred:
\begin{equation}
x'=r(t)\, x\,\left(1-\frac{x}{K(t)}\right)-\frac{a(t)\,x}{x+b(t)}\,,
\end{equation}
where the positively bounded-from-below maps $a$ and $b$ depend on the predator density and the average time between attacks of a predator.
In contrast to the Holling type II functional response term, which has strictly concave derivative on $\mathbb{R}$, the Holling type III functional response term $-a(t)\,x^2/(x^2+b(t))$ has strictly concave derivative on a smaller subset of $\mathbb{R}$. Some conditions to ensure strict d-concavity of the sum of a cubic Allee effect and the Holling type III term are in Ref.~\onlinecite{dno3}.

\section{Tipping and tracking}\label{sec:4tippingtracking}
In this section, we present the mathematical framework in which we describe critical transitions of the biological models of Section~\ref{sec:3biologicalmodels}.
In Subsection~\ref{subsec:parametershift}, we introduce the problem in a general setting.
Then, in Subsections~\ref{subsec:4Aconcavehypotheses} and \ref{subsec:4Bdconcavehypotheses}, we provide the particular concave and d-concave hypotheses and summarize all the dynamical possibilities of the transition equation \eqref{eq:x'=f(t,x,Gamma(t))}, which have been studied in previous papers.\cite{lno1,lno2,lno3,dno3}
Before that, we introduce the concepts of uniformly separated solutions, hyperbolic solutions and locally pullback attractive (resp. repulsive) solutions of nonautonomous equations \eqref{eq:2Ax'=h}, which will be repeatedly used.

\subsection{Concepts of nonautonomous attractivity}
Given an initial datum $(s,x_0)\in\mathbb{R}^2$, we denote by $t\mapsto x(t,s,x_0)$, the maximal solution of \eqref{eq:2Ax'=h} with initial condition $x_0$ at time $s\in\R$.
We say that $n$ bounded solutions $x_1,\,x_2,\,\dots,\,x_n$ of \eqref{eq:2Ax'=h}, with $n\geq2$, are \emph{uniformly separated} if $\inf_{t\in\mathbb{R}}|x_i(t)-x_j(t)|>0$ for any $1\leq i<j\leq n$.

A bounded solution $\tilde x\colon\mathbb{R}\rightarrow\mathbb{R}$ of \eqref{eq:2Ax'=h} is \emph{hyperbolic} if its variational equation $z'=h_x(t,\tilde x(t))\,z$ has an exponential dichotomy on $\mathbb{R}$, that is, if there exists $k\geq1$ and $\beta>0$ such that either
\begin{equation*}
\exp\int_s^t h_x(r,\tilde x(r))\, dr\leq ke^{-\beta(t-s)}\quad\text{for all }t\geq s
\end{equation*}
or
\begin{equation*}
\exp\int_s^t h_x(r,\tilde x(r))\, dr\leq ke^{\beta(t-s)}\quad\text{for all }t\leq s\,.
\end{equation*}
In the first case, we say that $\tilde x$ is \emph{attractive}, and, in the second one, that it is \emph{repulsive}.
The First Approximation Theorem ensures the uniform exponential attractivity at $+\infty$ or $-\infty$ of its graph.\cite{dno3}

Moreover, a solution $\tilde x\colon(-\infty,\beta)\rightarrow\mathbb{R}$ of \eqref{eq:2Ax'=h} is \emph{locally pullback attractive} if there exist $s_0<\beta$ and $\delta>0$ such that, if $s\leq s_0$ and $|x_0-\tilde x(s)|<\delta$, then $x(t,s,x_0)$ is defined for $t\in [s,s_0]$ and
\begin{equation*}
\lim_{s\rightarrow-\infty}\,\max_{x_0\in[\tilde x(s)-\delta,\tilde x(s)+\delta]}|\tilde x(t)-x(t,s,x_0)|=0\quad\text{for all }t\leq s_0\,.
\end{equation*}
In turn, a solution $\tilde x\colon(\alpha,\infty)\rightarrow\mathbb{R}$ of \eqref{eq:2Ax'=h} is \emph{locally pullback repulsive} if $\tilde y(t)=\tilde x(-t)$ is a locally pullback attractive solution of $y'=-h(-t,y)$.

\subsection{Parameter shifts and transition equations}\label{subsec:parametershift}
In addition to the intrinsic time variation of the system described in Section~\ref{sec:3biologicalmodels}, we assume that the population model under study undergoes a \emph{parameter shift} on an external parameter from one asymptotic value to another.\cite{alas,altw,awvc,lno1,lno2,lno3,dno3}
That is, we study the dynamics of a transition equation connecting two nonautonomous differential equations (quasiperiodic in the applications) which are referred to as the \emph{past} and the \emph{future equations}.
These past and future equations are given by some of the models described in Section~\ref{sec:3biologicalmodels}.

 Let $f\colon\mathbb{R}\times\mathbb{R}\times\mathbb{R}\rightarrow\mathbb{R}$, $(t,x,\gamma)\mapsto f(t,x,\gamma)$ be a continuous function whose regularity will be precised later on.
 The third variable represents the dependence on the external parameter $\gamma$.
 Two values $\gamma_-,\gamma_+$ of the external parameter define the past and future systems of our transition, that is, in the \emph{past system} the evolution of the state variable $x$ is governed by the nonautonomous differential equation \begin{equation}\label{eq:4pastequationgamma-}
     x'=f(t,x,\gamma_-)\,,
 \end{equation}
 while in the \emph{future system} it is governed by 
 \begin{equation}\label{eq:4futureequationgamma+}
     x'=f(t,x,\gamma_+)\,.
 \end{equation}
 
 Now, let us define how the system evolves from the past system to the future one. Let $\Gamma\colon\mathbb{R}\rightarrow\mathbb{R}$ be a function satisfying the following hypothesis, which will be always in force in both concave and d-concave cases,
 \begin{enumerate}[leftmargin=18pt]
\item[\hypertarget{CD}{\rm{\textbf{CD}}}] The map $\Gamma\colon\mathbb{R}\rightarrow\mathbb{R}$ is continuous and has asymptotic limits $\gamma_\pm=\lim_{t\rightarrow\pm\infty}\Gamma(t)$.
\end{enumerate}
The hypothesis \hyperlink{CD}{\rm{\textbf{CD}}} ensures that $\Gamma$ connects the past value $\gamma_-$ of the external parameter with the future one $\gamma_+$, giving rise to a \emph{transition equation}
 \begin{equation}\label{eq:4transitionequation}
 x'=f(t,x,\Gamma(t))
 \end{equation}
 whose right-hand side approaches $f(t,x,\gamma_+)$ for sufficiently large positive $t$ and $f(t,x,\gamma_-)$ for sufficiently large negative $t$.
 So, $\Gamma$ determines the parameter shift.
Throughout the paper, we will denote by $t\mapsto x_\Gamma(t,s,x_0)$ the maximal solution of \eqref{eq:4transitionequation} with initial condition $x_0$ at time $s\in\R$.
This equation describes the evolution of the number of individuals of the population through the parameter shift.
The dynamical possibilities of \eqref{eq:4transitionequation} are described in the following subsections.

\subsection{Concave framework}\label{subsec:4Aconcavehypotheses}
We say that the pair $(f,\Gamma)$ satisfies the family of properties \hypertarget{C}{\rm{\textbf{C}}} if $\Gamma$ satisfies \hyperlink{CD}{\rm{\textbf{CD}}} and $f$ satisfies the following properties:
\begin{enumerate}[leftmargin=18pt]
\item[\hypertarget{C1}{\rm{\textbf{C1}}}] $f$ is continuous, its derivative $f_x$ with respect to \emph{state variable} $x$ exists and is jointly continuous, and the restrictions $f,f_x\colon\mathbb{R}\times K_1\times K_2\rightarrow\mathbb{R}$ are bounded and uniformly continuous whenever $K_1$ and $K_2$ are compact subsets of $\mathbb{R}$.
\item[\hypertarget{C2}{\rm{\textbf{C2}}}] $f$ is \emph{coercive} in the sense that
\begin{equation*}
\limsup_{|x|\rightarrow\infty} \frac{f(t,x,\gamma)}{|x|}<0
\end{equation*}
uniformly for $(t,\gamma)\in\mathbb{R}\times \mathrm{cl}\,(\Gamma(\mathbb{R}))$.
\item[\hypertarget{C3}{\rm{\textbf{C3}}}] $f$ is strictly concave uniformly for $(t,\gamma)\in\mathbb{R}\times\mathrm{cl}\,(\Gamma(\mathbb{R}))$, that is, for any $j\in\mathbb{N}$, there exists $\delta_j>0$ such that
\begin{equation}\label{eq:2Dstrictconcavity}
    \sup_{x_1,\,x_2\in[-j,j]\,,\;x_1\neq x_2} \frac{f_{x}(t,x_2,\gamma)-f_{x}(t,x_1,\gamma)}{x_2-x_1}<-\delta_j
\end{equation}
for all $(t,\gamma)\in\mathbb{R}\times\mathrm{cl}\,(\Gamma(\mathbb{R}))$.
\item[\hypertarget{C4}{\rm{\textbf{C4}}}] $x'=f(t,x,\gamma)$ has two uniformly separated hyperbolic solutions for $\gamma\in\{\gamma_-,\gamma_+\}$, which we will call $\tilde r_{\gamma_-}<\tilde a_{\gamma_-}$ and $\tilde r_{\gamma_+}<\tilde a_{\gamma_+}$ respectively.
\end{enumerate}

Hypothesis \hyperlink{C1}{\rm{\textbf{C1}}} establishes the type of time dependency allowed (required to use the skewproduct formalism of Subsection~\ref{subsec:2Bskewproductpreliminaries}).
We remark that, although we are working within the continuous framework for the sake of simplicity, most of the results of the present paper can be stated in Carathéodory topologies,\cite{lno1carat,lno2carat,lno3carat} as in Refs.~\onlinecite{lno2,lno3}.
The simplest examples of strictly concave functions are quadratic polynomials with a negatively bounded-from-above quadratic coefficient.

When conditions \hyperlink{CD}{{\rm\textbf{CD}}}, \hyperlink{C1}{{\rm\textbf{C1}}} and \hyperlink{C2}{{\rm\textbf{C2}}} hold, there exist a solution $r_\Gamma$ delimiting from below the set of
solutions which are bounded as time increases (if there is any), and 
a solution $a_\Gamma$ delimiting from above the set of solutions which 
are bounded as time decreases (if there is any). In fact, one of these solutions is 
globally defined if and only if the set 
$\mathcal B$ of bounded solutions of \eqref{eq:4transitionequation} is nonempty, in which case also the other one is 
globally defined and
\begin{equation*}
\mathcal{B}=\{(s,x_0)\colon\; r_\Gamma(s)\leq x_0\leq a_\Gamma(s)\}\,.
\end{equation*}
These properties are proved in Theorem 3.1 of Ref.~\onlinecite{lno3}.
On the other hand, if, in addition, condition \hyperlink{C3}{{\rm\textbf{C3}}} holds, then Theorem 3.4 of Ref.~\onlinecite{lno3} ensures that, for \eqref{eq:4transitionequation} as well as for the past and future equations \eqref{eq:4pastequationgamma-} and \eqref{eq:4futureequationgamma+}, the existence of two uniformly 
separated solutions is equivalent to the existence of two hyperbolic solutions, 
in which case the hyperbolic solutions are the uniformly separated ones, the upper one is attractive, 
and the lower one repulsive.
This is the 
situation for the hyperbolic solutions $\tilde r_{\gamma_\pm}<\tilde a_{\gamma_\pm}$ if, in addition, 
\hyperlink{C4}{{\rm\textbf{C4}}} holds. 
%An analogous proof to that of Proposition 3.16 of Ref.~\onlinecite{dno3} ensures that, if \hyperlink{C3}{\rm{\textbf{C3}}} and \hyperlink{C4}{\rm{\textbf{C4}}} are satisfied, then the right-hand side of \eqref{eq:4transitionequation} satisfies the strict concavity hypothesis of Subsection~\ref{subsec:2Dstrictconcavityanddconcavity}, so the conclusions of Theorem~3.5 of Ref.~\onlinecite{lno1} can be drawn.
%Analogous arguments to those of the proof of Theorem 3.3 of Ref.~\onlinecite{dno3} show that the hyperbolic solutions $\tilde r_{\gamma_-}<\tilde a_{\gamma_-}$ and $\tilde r_{\gamma_+}<\tilde a_{\gamma_+}$ given by \hyperlink{C5}{\rm{\textbf{C5}}} are uniformly separated, and that

The family of hypotheses \hyperlink{C}{{\rm\textbf{C}}} allows us to repeat the arguments of Theorem 4.4 of Ref.~\onlinecite{lno3}
in order to check the statements of the rest of the subsection.
In particular, $\lim_{t\to-\infty}(a_\Gamma(t)-\tilde a_{\gamma_-}(t))=0$, $a_\Gamma$ is locally pullback attractive,
$\lim_{t\to\infty}(r_\Gamma(t)-\tilde r_{\gamma_+}(t))=0$, and $r_\Gamma$ is locally pullback repulsive.
The properties of $a_\Gamma$ and $r_\Gamma$ 
determine the three different possibilities for the global dynamics of \eqref{eq:4transitionequation} under the family of hypotheses \hyperlink{C}{{\rm\textbf{C}}}:
\begin{itemize}[leftmargin=*]
\item[-] It is in \textsc{Case} \hypertarget{Ac}{{\rm A}} if there are two different hyperbolic solutions, in which case they are $\tilde a_\Gamma=a_\Gamma$ and $\tilde r_\Gamma=r_\Gamma$. 
In addition, $\tilde a_\Gamma$ is attractive and satisfies $\lim_{t\rightarrow\infty}(\tilde a_\Gamma(t)-\tilde a_{\gamma_+}(t))=0$, and $\tilde r_\Gamma$ is repulsive and satisfies $\lim_{t\rightarrow-\infty}(\tilde r_\Gamma(t)-\tilde r_{\gamma_-}(t))=0$.
\item[-] It is \textsc{Case} \hypertarget{Bc}{{\rm B}} if there is exactly one bounded (and hence nonhyperbolic) solution, in which case it is given by $r_\Gamma=a_\Gamma$, and hence it is locally pullback attractive and locally pullback repulsive.
\item[-] It is in \textsc{Case} \hypertarget{Cc}{{\rm C}} if there are no bounded solutions.
\end{itemize}
%We refer to \textsc{Case} \hyperlink{Ac}{{\rm A}} as \emph{tracking} and to \textsc{Cases} \hyperlink{Bc}{{\rm B}} and \hyperlink{Cc}{{\rm C}} as \emph{tipping}. 
Their respective dynamics are identical to those described in detail in Theorems~4.4, 4.5 and 4.6 of Ref.~\onlinecite{lno1} respectively.
In particular, \textsc{Case} \hyperlink{Ac}{{\rm A}} holds (resp. \textsc{Case} \hyperlink{Bc}{{\rm B}}) if and only if
$a_\Gamma(t_0)>r_\Gamma(t_0)$ (resp. $a_\Gamma(t_0)=r_\Gamma(t_0)$) for a 
$t_0\in \mathbb R$, in which case the inequality (resp. equality) holds for all $t\in\mathbb R$.

In all the cases, the backward and forward behavior 
of any solution $\bar x_\Gamma(t)$ of \eqref{eq:4transitionequation} is determined by its relation with $r_\Gamma$ and $a_\Gamma$. 
As said before, $\bar x_\Gamma(t)$ is bounded as time increases if and only if there exists $s$ in the definition domain of $r_\Gamma$ such that $\bar x_\Gamma(s)\ge r_\Gamma(s)$. 
In addition, if the inequality is strict, then $\lim_{t\to\infty}(\bar x_\Gamma(t)-\tilde a_{\gamma_+}(t))=0$.
Analogously, $\bar x_\Gamma$ is bounded as time decreases if and only if there exists $s$ in the definition domain of $a_\Gamma$ such that $\bar x_\Gamma(s)\le a_\Gamma(s)$, and, if the inequality is strict, then $\lim_{t\to-\infty}(x_\Gamma(t)- \tilde r_{\gamma_-}(t))=0$.
\subsection{D-concave framework}\label{subsec:4Bdconcavehypotheses}
We say that the pair $(f,\Gamma)$ satisfies the family of properties \hypertarget{D}{\rm{\textbf{D}}} if $\Gamma$ satisfies \hyperlink{CD}{{\rm\textbf{CD}}} and $f$ satisfies the following properties:
\begin{enumerate}[leftmargin=18pt]
\item[\hypertarget{D1}{\rm{\textbf{D1}}}] $f$ is continuous, its derivatives $f_x$ and $f_{xx}$ with respect to \emph{state variable} $x$ exist and are jointly continuous, and the restrictions $f,f_x,f_{xx}\colon\mathbb{R}\times K_1\times K_2\rightarrow\mathbb{R}$ are bounded and uniformly continuous whenever $K_1$ and $K_2$ are compact subsets of $\mathbb{R}$.
\item[\hypertarget{D2}{\rm{\textbf{D2}}}] $f$ is \emph{coercive}, in the sense that
\begin{equation*}
\lim_{|x|\rightarrow\infty} \frac{f(t,x,\gamma)}{x}=-\infty
\end{equation*}
uniformly for $(t,\gamma)\in\mathbb{R}\times \mathrm{cl}\,(\Gamma(\mathbb{R}))$.
\item[\hypertarget{D3}{\rm{\textbf{D3}}}] $f$ is strictly d-concave uniformly for $(t,\gamma)\in\mathbb{R}\times\mathrm{cl}\,(\Gamma(\mathbb{R}))$, that is, for any $j\in\mathbb{N}$, there exists $\delta_j>0$ such that
\begin{equation}\label{eq:2Dstrictdconcavity}
    \sup_{x_1,\,x_2\in[-j,j]\,,\;x_1\neq x_2} \frac{f_{xx}(t,x_2,\gamma)-f_{xx}(t,x_1,\gamma)}{x_2-x_1}<-\delta_j
\end{equation}
for all $(t,\gamma)\in\mathbb{R}\times\mathrm{cl}\,(\Gamma(\mathbb{R}))$.
\item[\hypertarget{D4}{\rm{\textbf{D4}}}] $x'=f(t,x,\gamma)$ has three different uniformly separated hyperbolic solutions for $\gamma\in\{\gamma_-,\gamma_+\}$,
which we will call $\tilde l_{\gamma_-}<\tilde m_{\gamma_-}<\tilde u_{\gamma_-}$ and $\tilde l_{\gamma^+}<\tilde m_{\gamma_+}<\tilde u_{\gamma_+}$ respectively.
\end{enumerate}

The simplest examples of strictly d-concave functions are cubic polynomials with a negatively bounded-from-above cubic coefficient, although all the d-concave numerical experiments of this paper are done in the nonpolynomial case.

Analogously to the concave framework, if $(f,\Gamma)$ satisfies the family of hypothesis \hyperlink{D}{{\rm{\textbf{D}}}}, then Theorem 3.7 of Ref.~\onlinecite{dno3} ensures the existence of a unique solution $l_\Gamma$ of \eqref{eq:4transitionequation} satisfying $\lim_{t\rightarrow-\infty}(l_\Gamma(t)-\tilde l_{\gamma_-}(t))=0$, a unique solution $m_\Gamma$ satisfying $\lim_{t\rightarrow\infty}(m_\Gamma(t)-\tilde m_{\gamma_+}(t))=0$, and a unique solution $u_\Gamma$ satisfying $\lim_{t\rightarrow-\infty}(u_\Gamma(t)-\tilde u_{\gamma_-}(t))=0$.
The solutions $l_\Gamma$ and $u_\Gamma$ are locally pullback attractive and globally defined while $m_\Gamma$ is locally pullback repulsive and not necessarily globally defined.
We call $u_\Gamma$ (resp. $l_\Gamma$) the upper (resp. lower) locally pullback attractive solution of \eqref{eq:4transitionequation}.
We will respectively denote them by $\tilde l_\Gamma$, $\tilde m_\Gamma$ and $\tilde u_\Gamma$ if they are hyperbolic.
Moreover, the hyperbolic solutions $\tilde l_{\gamma_-},\tilde u_{\gamma_-},\tilde l_{\gamma_+},\tilde u_{\gamma_+}$ given by \hyperlink{D4}{{\rm\textbf{D4}}} are attractive, while $\tilde m_{\gamma_-},\tilde m_{\gamma_+}$ are repulsive.

Furthermore, all the solutions of the transition equation \eqref{eq:4transitionequation} are globally forward defined, and the set $\mathcal{B}$ of bounded solutions of \eqref{eq:4transitionequation} is always nonempty and takes the form:
\begin{equation*}
    \mathcal{B}=\{(s,x_0)\colon\;l_\Gamma(s)\leq x_0\leq u_\Gamma(s)\}\,.
\end{equation*}

The dynamics of \eqref{eq:4transitionequation} fits one of the following scenarios:\cite{dno3}
\begin{itemize}[leftmargin=*]
\item \textsc{Case} \hypertarget{Ad}{{\rm A}} if there are three hyperbolic solutions, which are the unique three
uniformly separated solutions and which are precisely $\tilde l_\Gamma<\tilde m_\Gamma<\tilde u_\Gamma$. Then, $\tilde l_\Gamma$ is hyperbolic attractive and $\lim_{t\rightarrow\infty}(\tilde l_\Gamma(t)-\tilde l_{\gamma_+}(t))=0$, $\tilde m_\Gamma$ is hyperbolic repulsive and $\lim_{t\rightarrow-\infty}(\tilde m_\Gamma(t)-\tilde m_{\gamma_-}(t))=0$, and $\tilde u_\Gamma$ is hyperbolic attractive and $\lim_{t\rightarrow\infty}(\tilde u_\Gamma(t)-\tilde u_{\gamma_+}(t))=0$.
\item \textsc{Case} \hypertarget{Bd}{{\rm B}} if there are exactly two uniformly separated solutions,
one of them being the only hyperbolic solution, of attractive type, and the other one being nonhyperbolic, locally pullback attractive and locally pullback repulsive.
There exist two symmetric subcases \hypertarget{B1d}{{\rm B1}} and \hypertarget{B2d}{{\rm B2}}.
In \textsc{Case} \hyperlink{B1d}{{\rm B1}} (resp. \hyperlink{B2d}{{\rm B2}}), $\tilde u_\Gamma$ (resp. $\tilde l_\Gamma$) is the hyperbolic solution and $m_\Gamma=l_\Gamma$ (resp. $m_\Gamma=u_\Gamma$) is the nonhyperbolic one; and $\lim_{t\rightarrow\infty}(\tilde u_\Gamma(t)-\tilde u_{\gamma_+}(t))=0$ (resp. $\lim_{t\rightarrow\infty}(\tilde l_\Gamma(t)-\tilde l_{\gamma_+}(t))=0$).
\item \textsc{Case} \hypertarget{Cd}{{\rm C}} if there are no uniformly separated solutions and there are exactly two hyperbolic solutions,
which are attractive and correspond to $\tilde l_\Gamma$ and $\tilde u_\Gamma$.
There exist two symmetric subcases \hypertarget{C1d}{{\rm C1}} and \hypertarget{C2d}{{\rm C2}}.
In \textsc{Case} \hyperlink{C1d}{{\rm C1}},
$\lim_{t\rightarrow\infty}(\tilde u_\Gamma(t)-\tilde u_{\gamma_+}(t))=\lim_{t\rightarrow\infty}(\tilde l_\Gamma(t)-\tilde u_{\gamma_+}(t))=0$,
while in \textsc{Case} \hyperlink{C2d}{{\rm C2}},
$\lim_{t\rightarrow\infty}(\tilde u_\Gamma(t)-\tilde l_{\gamma_+}(t))=\lim_{t\rightarrow\infty}(\tilde l_\Gamma(t)-\tilde l_{\gamma_+}(t))=0$. In both subcases $m_\Gamma$ is unbounded.
\end{itemize}
Theorems~3.7, 3.9 and 3.10 of Ref.~\onlinecite{dno3} accurately describe the three cases.
In all the cases, $m_\Gamma$ governs the forward behavior of solutions, that is, given $s\in\mathbb{R}$ in the interval of definition of $m_\Gamma$ (recall that $m_\Gamma$ is globally defined in \textsc{Cases} \hyperlink{Ad}{{\rm A}} and \hyperlink{Bd}{{\rm B}}), $\lim_{t\rightarrow\infty} (x_\Gamma(t,s,x_0)-\tilde u_{\gamma_+}(t))=0$ if $x_0>m_\Gamma(s)$, and $\lim_{t\rightarrow\infty} (x_\Gamma(t,s,x_0)-\tilde l_{\gamma_+}(t))=0$ if $x_0<m_\Gamma(s)$.
In addition, for any $s\in\mathbb{R}$, $\lim_{t\rightarrow-\infty}(x_\Gamma(t,s,x_0)-\tilde m_{\gamma_-}(t))=0$ if and only if $x_0\in(l_\Gamma(s),u_\Gamma(s))$.

\subsection{Mechanisms producing critical transitions}\label{subsec:4Cmechanisms}
Distinct physical mechanisms have been identified as sources of tipping in applied sciences.\cite{altw,aajqw,awvc}
Now, we introduce some parametric families of parameter shifts $c\mapsto\Gamma^c$, with $\Gamma^c$ satisfying \hyperlink{CD}{{\rm\textbf{CD}}} for all the considered values of $c$, which represent the influence of these different mechanisms on the system, and we define what a critical transition is in our models.
We will study the parametric family of transition equations
\begin{equation}\label{eq:4Cparametrictransitionequation}
x'=f(t,x,\Gamma^c(t))\,.
\end{equation}
Hereafter, we will denote the equation of the parametric family \eqref{eq:4Cparametrictransitionequation} for a fixed value $c$ of the parameter by \eqref{eq:4Cparametrictransitionequation}$_c$, the locally pullback repulsive and attractive solutions of \eqref{eq:4Cparametrictransitionequation}$_c$ under the family of hypotheses \hyperlink{C}{{\rm\textbf{C}}} by $r_c$ and $a_c$, the locally pullback repulsive and attractive solutions of \eqref{eq:4Cparametrictransitionequation}$_c$ under the family of hypotheses \hyperlink{D}{{\rm\textbf{D}}} by $l_c$, $m_c$ and $u_c$,
and the maximal solution of \eqref{eq:4Cparametrictransitionequation}$_c$ by $t\mapsto x_c(t,s,x_0)$ with initial condition $x_0$ at time $s\in\R$

The considered mechanisms are the following:
\begin{itemize}[leftmargin=*]
\item \emph{Rate parametric problem,} $\Gamma^c(t)=\Gamma(ct)$: we study %the parametric family
\begin{equation}\label{eq:4Cparametrictransitionrate}
    x'=f(t,x,\Gamma(ct))\,
\end{equation}
so that $c>0$ represents the rate of the transition,
as it modifies the speed at which the parameter shift takes place.
\item \emph{Phase parametric problem,} $\Gamma^c(t)=\Gamma(t+c)$: we study %the parametric family 
\begin{equation}\label{eq:4Cparametrictransitionphase}
x'=f(t,x,\Gamma(t+c))\,
\end{equation}
so that different values of $c\in\mathbb{R}$ mean different transition phases,
as it is a delay or advance in time of the parameter shift.
\item \emph{Size parametric problem,} $\Gamma^c(t)=c\,\Gamma(t)$: we study %the parametric family
\begin{equation}\label{eq:4Cparametrictransitionsize}
x'=f(t,x,c\,\Gamma(t))\,
\end{equation}
restricting ourselves to $f(t,x,\gamma)=h(t,x-\gamma)$. In this case, the future equation is a lift of the past one, and $c\in\mathbb{R}$ determines the size of the transition, that is,
how far are the maximum and minimum values reached by $\Gamma^c$.
\end{itemize}
\begin{definition}\label{def:3.1tippingtracking}
Let $(f,\Gamma)$
satisfy the family of hypotheses \hyperlink{C}{{\rm\textbf{C}}} or \hyperlink{D}{{\rm\textbf{D}}}.
We shall say that
\eqref{eq:4transitionequation} shows {\em tracking} if it is in \hyperlink{Ac}{{\sc Case A}} and that \eqref{eq:4transitionequation} shows {\em tipping} (or equivalently undergoes a {\em critical transition}) if it is either in \hyperlink{Bc}{{\sc Cases B}} or \hyperlink{Cc}{{\sc C}}.

    We shall say that the parametric family \eqref{eq:4Cparametrictransitionequation} undergoes a \emph{rate-induced, phase-induced} or \emph{size-induced critical transition} if there exist $c_1\neq c_2$ such that \eqref{eq:4Cparametrictransitionequation}$_{c_1}$ and \eqref{eq:4Cparametrictransitionequation}$_{c_2}$ are in different dynamical cases.

    We shall say $c_0$ is a {\em tipping point} for \eqref{eq:4Cparametrictransitionequation} if the parametric family changes of dynamical case at this value.
\end{definition}
%Special attention will be paid to \emph{transversal} critical transitions, that is, cases in which there exists a value of the parameter of \textsc{Case} \hyperlink{Bc}{{\rm B}} which separates an open region of \textsc{Case} \hyperlink{Ac}{{\rm A}} from an open region of \textsc{Case} \hyperlink{Cc}{{\rm C}}.
This notion of tracking corresponds somehow to the notion of end-point tracking in Ref. \onlinecite{aspw}, that is, the local attractors in the past perfectly connect with the local attractors in the future.
%In both {\color{red} concave} and d-concave frameworks, \textsc{Cases} \hyperlink{Ac}{{\rm A}} and \hyperlink{Cc}{{\rm C}} are robust situations which persist under sufficiently small perturbations of the parameter shift, while 
\textsc{Case} \hyperlink{Bc}{{\rm B}} is frequently a highly unstable situation which separates \textsc{Cases} \hyperlink{Ac}{{\rm A}} and \hyperlink{Cc}{{\rm C}}, meaning that the dynamical case of \eqref{eq:4Cparametrictransitionequation}$_c$ is not robust to perturbation in $c$ at those points.
In this situation, we will say that a \emph{transversal} critical transition takes place.
Thus, the change between \textsc{Case} \hyperlink{Ac}{{\rm A}} and \hyperlink{Cc}{{\rm C}} as the parameter $c$ varies in \eqref{eq:4Cparametrictransitionequation} is triggered by a nonautonomous saddle-node bifurcation: as the parameter increases or decreases, two hyperbolic solutions approach, collide at the tipping point $c_0$ and disappear.
This question will be studied in the following subsection in a different framework.

\subsection{Concave and d-concave bifurcations}\label{subsec:4Dbifurcations}
In the biological models considered in this paper, a saddle-node $\gamma$-parametric bifurcation of the intermediate equations 
\begin{equation}\label{eq:4intermediate}
    x'=f(t,x,\gamma)
\end{equation}
can rigorously justify a critical transition of \eqref{eq:4transitionequation}.
In this subsection, we present rigorous theorems describing how this saddle-node bifurcation takes place
for the particular monotone $\gamma$-parametric variation,
\begin{equation}\label{eq:4Dlambdaparamteric}
    x'=h(t,x)+\gamma\,.
\end{equation}
Nevertheless, we remark that the result we show serve as a pattern for more general problems.
We emphasize the differences between the study of parametric problems \eqref{eq:4Cparametrictransitionequation} and \eqref{eq:4Dlambdaparamteric}: problem \eqref{eq:4Cparametrictransitionequation} is more general than \eqref{eq:4Dlambdaparamteric} in the sense that the parametric variation is not necessarily monotonic. However, \eqref{eq:4Dlambdaparamteric} is more general than \eqref{eq:4Cparametrictransitionequation} in the sense that it does not necessarily correspond to a transition equation with known asymptotic dynamics.
%We remark that the parametric problem \eqref{eq:4Cparametrictransitionequation} is more general than \eqref{eq:4Dlambdaparamteric} in the sense that the parametric variation is not necesarily monotonic, and less general in the sense

%which serves as a standard to understand the nonautonomous bifurcations appearing in the parametric variations described in Subsection~\ref{subsec:4Cmechanisms} (which are more general in the sense that the right-hand side of \eqref{eq:4Cparametrictransitionequation} is not monotonous on $c$, but less general in the sense that they are a parameter shift between to systems with known dynamics while \eqref{eq:4Dlambdaparamteric} allows for a more complicated structure in the skewproduct base).
In what follows, $C^{0,i}$, for $i\in\{1,2\}$, stands for the space of continuous functions $h\colon\mathbb{R}\times\mathbb{R}\rightarrow\mathbb{R}$, $(t,x)\mapsto h(t,x)$ such that its derivatives up to order $i$ with respect to the second variable exist (denoted by $h_x$ and $h_{xx}$), are continuous, and their restrictions to $\mathbb{R}\times K$ (also that of $h$) are bounded and uniformly continuous for any compact $K\subset\mathbb{R}$.

The following theorem, which describes the bifurcation in a concave setting, is proved in Theorem~3.5 of Ref.~\onlinecite{lno3}.
\begin{theorem}[Concave saddle-node bifurcation]\label{th:4Dconcavebifurcationtheorem} Let $h\in C^{0,1}$ be strictly concave uniformly for $t\in\mathbb{R}$ (see condition \hyperlink{C3}{{\rm\textbf{C3}}}) and assume that there exists $\delta>0$ such that $\limsup_{|x|\rightarrow\infty} h(t,x)/|x|<0$ uniformly for $t\in\mathbb{R}$.
Under these assumptions, whenever the set $\mathcal{B}_\gamma$ of bounded solutions of \eqref{eq:4Dlambdaparamteric}$_\gamma$ is nonempty, it is bounded and there exist two solutions $r_\gamma\leq a_\gamma$ such that
\begin{equation*}
 \mathcal{B}_\gamma=\{(s,x_0)\colon\, r_\gamma(s)\leq x_0\leq a_\gamma(s)\}\,.
\end{equation*}
Then, there exists $\gamma^*\in\mathbb{R}$ such that
\begin{enumerate}[label=\rm{(\roman*)},leftmargin=18pt]
\item $\mathcal{B}_\gamma$ is empty if and only if $\gamma<\gamma^*$.
\item For $\gamma=\gamma^*$, $\inf_{t\in\mathbb{R}}(a_{\gamma^*}(t)-r_{\gamma^*}(t))=0$ and there are no hyperbolic solutions.
\item For $\gamma>\gamma^*$, there exist only two hyperbolic solutions $\tilde r_\gamma<\tilde a_\gamma$ of \eqref{eq:4Dlambdaparamteric}$_\gamma$, which are uniformly separated. Moreover, $\tilde r_\gamma= r_\gamma$ is repulsive, $\tilde a_\gamma=a_\gamma$ is attractive, and for any fixed $t\in\mathbb{R}$, $\gamma\mapsto\tilde a_\gamma(t)$ is a strictly increasing map on $(\gamma^*,\infty)$ and $\gamma\mapsto\tilde r_\gamma(t)$ is a strictly decreasing map on $(\gamma^*,\infty)$, with $\lim_{\gamma\downarrow\gamma^*}\tilde a_\gamma(t)=a_{\gamma^*}(t)$ and $\lim_{\gamma\downarrow\gamma^*}\tilde r_\gamma(t)=r_{\gamma^*}(t)$.
\end{enumerate}
\end{theorem}
In short, Theorem \ref{th:4Dconcavebifurcationtheorem} describes the existence of two uniformly separated hyperbolic solutions, attractive the upper one and repulsive the lower one, which approach as $\gamma\downarrow\gamma_*$. At $\gamma_*$ they are not uniformly separated (but not necessarily collide at any point of the real line), and they dissapear for $\gamma<\gamma_*$.

The following theorem describes the saddle-node nonautonomous bifurcation of \eqref{eq:4Dlambdaparamteric} in a d-concave setting. Most of its proof relies on the proving arguments of Theorem 5.10 of Ref.~\onlinecite{dno1}, although some complementary details, which are necessary in this more general case, are given in Appendix \ref{appendix2}.
\begin{theorem}[D-concave saddle-node bifurcations]\label{th:4Dd-concavebifurcationtheorem} Let $h\in C^{0,2}$ be strictly d-concave uniformly for $t\in\mathbb{R}$ (see condition \hyperlink{D3}{{\rm\textbf{D3}}}) and assume that $\lim_{|x|\rightarrow\infty} h(t,x)/x=-\infty$ uniformly for $t\in\mathbb{R}$.
Under these assumptions, the set $\mathcal{B}_\gamma$ of bounded solutions of \eqref{eq:4Dlambdaparamteric}$_\gamma$ is nonempty and bounded, and there exist two solutions $l_\gamma\leq u_\gamma$ such that
\begin{equation*}
    \mathcal{B}_\gamma=\{(s,x_0)\colon\; l_\gamma(s)\leq x_0\leq u_\gamma(s)\}\,.
\end{equation*}
Assume that there exist $\gamma_0\in\mathbb{R}$ such that \eqref{eq:4Dlambdaparamteric}$_{\gamma_0}$ has three uniformly separated solutions.
Then, there exists an interval $I=(\gamma_1,\gamma_2)$ with $\gamma_0\in I$ such that
\begin{enumerate}[label=\rm{(\roman*)},leftmargin=18pt]
\item for $\gamma\in I$, there exist only three uniformly separated solutions $\tilde l_\gamma<\tilde m_\gamma<\tilde u_\gamma$. In addition, $\tilde l_\gamma=l_\gamma$ and $\tilde u_\gamma=u_\gamma$ are hyperbolic attractive, $\tilde m_\gamma$ is hyperbolic repulsive, and for any fixed $t\in\mathbb{R}$, $\gamma\mapsto \tilde l_\gamma(t),\tilde u_\gamma(t)$ are strictly increasing maps on $I$ and $\gamma\mapsto \tilde m_\gamma(t)$ is strictly decreasing on $I$.
\item For $\gamma=\gamma_2$ (resp. $\gamma=\gamma_1$), the solution $\tilde u_{\gamma_2}=u_{\gamma_2}$ (resp. $\tilde l_{\gamma_1}=l_{\gamma_1}$) is hyperbolic attractive, and $\inf_{t\in\mathbb{R}}(m_{\gamma_2}(t)-l_{\gamma_2}(t))=0$ (resp. $\inf_{t\in\mathbb{R}}(u_{\gamma_1}(t)-m_{\gamma_1}(t))=0$), where $m_{\gamma_2}(t)=\lim_{\gamma\uparrow\gamma_2}\tilde m_\gamma(t)$ (resp. $\tilde m_{\gamma_1}(t)=\lim_{\gamma\downarrow\gamma_1}\tilde m_\gamma(t)$).
\item For $\gamma>\gamma_2$ (resp. $\gamma<\gamma_1$), $\tilde u_\gamma=u_\gamma$ (resp. $\tilde l_\gamma=l_\gamma$) is hyperbolic attractive and $\lim_{t\rightarrow\infty}(\tilde u_\gamma(t)-l_\gamma(t))=0$ (resp. $\lim_{t\rightarrow\infty}(u_\gamma(t)-\tilde l_\gamma(t))=0$), that is, there are no uniformly separated solutions.
\end{enumerate}
\end{theorem}
To sum up, the previous theorem describes a two-sided phenomenon of saddle-node bifurcations. Starting from the value $\gamma_0$ of the parameter for which we have three uniformly separated hyperbolic solutions, the two upper hyperbolic solutions approach and undergo a saddle-node bifurcation as $\gamma\uparrow\gamma_2$, and the two lower ones approach and undergo another saddle-node bifurcation as $\gamma\downarrow\gamma_1$.

The nonautonomous equation \eqref{eq:4Dlambdaparamteric} can be studied through topological dynamics via the construction of a continuous skew-product flow on the hull $\Omega_h$ of $h$ (see Subsection~\ref{subsec:2Bskewproductpreliminaries} in Appendix~\ref{appendix1}). 
Assuming that $h$ is recurrent, i.e., that every orbit on $\Omega_h$ is dense, Theorems \ref{th:4Dconcavebifurcationtheorem} and \ref{th:4Dd-concavebifurcationtheorem} can be applied to any equation \eqref{eq:2Bfamilyskewproduct}$_\omega$ for $\omega\in\Omega_h$.
In this case, it is known that the two semicontinuous equilibria at the saddle-node bifurcation points which extend to $\Omega_h$ any of the pairs $(r_{\gamma^*},a_{\gamma^*})$, $(l_{\gamma_2},m_{\gamma_2})$ or $(m_{\gamma_1},u_{\gamma_1})$ (depending on which case we are)
are equal on the residual subset of $\Omega_h$ of common continuity points.\cite{no1,dno1}

\section{Finite-time Lyapunov exponents and early-warning signals}\label{sec:5earlywarningsignals}
In this section, we present a construction of a skewproduct base $\Omega$ on which we can jointly represent all the equations of the parametric family \eqref{eq:4Cparametrictransitionequation}, where $\Omega$ can be understood as a heteroclinic orbit connecting compact invariant sets representing past $\Omega_-$ and future $\Omega_+$.
We also show that finite-time Lyapunov exponents are effective early-warning signals for some critical transitions (including the transversal ones) occurring for one-parametric families of equations \eqref{eq:4Cparametrictransitionequation} which satisfy the families of hypotheses \hyperlink{C}{{\rm\textbf{C}}} or \hyperlink{D}{{\rm\textbf{D}}}, and for which the parameter $c$ means either the rate, the phase or the size of the parameter shift.

\subsection{Joint skewproduct structure}\label{subsec:4jointskewproduct}
In this subsection, we present the joint skewproduct formulation for the parametric problems \eqref{eq:4Cparametrictransitionequation} presented in Subsection~\ref{subsec:4Cmechanisms}.
The definition and construction of the skewproduct formalism from a single vector field can be found in Subsection~\ref{subsec:2Bskewproductpreliminaries} of Appendix~\ref{appendix1}.

For the rate and phase mechanisms, we get parametric families \eqref{eq:4Cparametrictransitionequation} that have the same past and future equations for all the values of the parameter. 
As long as there is some $c_0$ for which \eqref{eq:4Cparametrictransitionequation}$_{c_0}$ is in \textsc{Case} \hyperlink{Bc}{{\rm B}}, we are able to construct a joint skewproduct representation of all the equations of the parametric family \eqref{eq:4Cparametrictransitionequation} on a common base $\Omega$, which is defined as the hull of $f_{\Gamma^{c_0}}(t,x)=f(t,x,\Gamma^{c_0}(t))$.

In \textsc{Case} \hyperlink{Bc}{{\rm B}}, the dynamics of the transition equation is essentially different from that of the past and future equations (e.g. it has a different number of hyperbolic solutions), so $\{(f_{\Gamma^{c_0}}){\cdot}s\colon s\in\mathbb R\}$ is an orbit that is neither contained in its $\boldsymbol\alpha$-limit set $\Omega_-$ (past) nor in its $\boldsymbol\omega$-limit set $\Omega_+$ (future): it can be understood as a heteroclinic orbit connecting these two compact invariant sets.
In fact, as provided in Proposition~\ref{prop:AstructureofOmegatransitive},
\begin{equation*}
    \Omega=\Omega_-\cup\big\{(f_{\Gamma^{c_0}}){\cdot}s\colon s\in\mathbb R\big\}\cup\Omega_+\,.
\end{equation*}
Then, we can use $\Omega$ to represent all the $c$-parametric family of transition equations \eqref{eq:4Cparametrictransitionequation}.
Specifically, denoting $f_{\Gamma^c}(t,x)=f(t,x,\Gamma^c(t))$, if \eqref{eq:4Cparametrictransitionequation} is a rate \eqref{eq:4Cparametrictransitionrate} or phase \eqref{eq:4Cparametrictransitionphase} problem, then we define
\begin{equation*}\label{eq:5Askewfunctionratephase}
 \mathfrak{f}_{\Gamma^c}(\omega,x)=\left\{\begin{array}{ll}
f_{\Gamma^c}(s,x)\quad&\text{if }\omega=f_{\Gamma^{c_0}}{\cdot}s\text{ for some }s\in\mathbb{R}\,,\\[1ex]
\omega(0,x)&\text{otherwise; i.e., if $\omega\in\Omega_\pm$}\,,
\end{array}\right.
\end{equation*}
check its continuity, and consider the $c$-parametric family of families
\begin{equation}\label{eq:5Atransitionskewproduct}
 x'=\mathfrak{f}_{\Gamma^c}(\omega{\cdot}t,x)\,,\quad\omega\in\Omega\,. 
\end{equation}
So, \eqref{eq:4Cparametrictransitionequation} corresponds to the previous equation for $\omega=f_{\Gamma^{c_0}}$.
We remark that $\mathfrak f_{\Gamma^c}(\omega,x)$ does not depend on $c$ for $\omega\in\Omega_\pm$: its definition only changes with $c$ on the orbit of $f_{\Gamma^{c_0}}$.

For the size case \eqref{eq:4Cparametrictransitionsize}, the past and future equations take the form 
\begin{equation}
    x'=f(t,x-c\,\gamma_\pm)\,,
\end{equation}
so the corresponding dynamics for all $c$ are lifts of that of $x'=f(t,x)$. 
A similar construction can be made in this case, but we omit it because it is a bit more technical.

This common framework elucidates the topological structure underlying the transition equation, makes clear that the Lyapunov exponents (see Appendix \ref{subsec:2CLyapunovexponents}) of past and future equations are unaffected by the parametric variation, and motivates the calculation of the finite-time Lyapunov exponents on the transition equation.

\subsection{Finite-time Lyapunov exponents}\label{subsec:5BFTLE}
\begin{figure*}
\centering
\includegraphics[width=\textwidth]{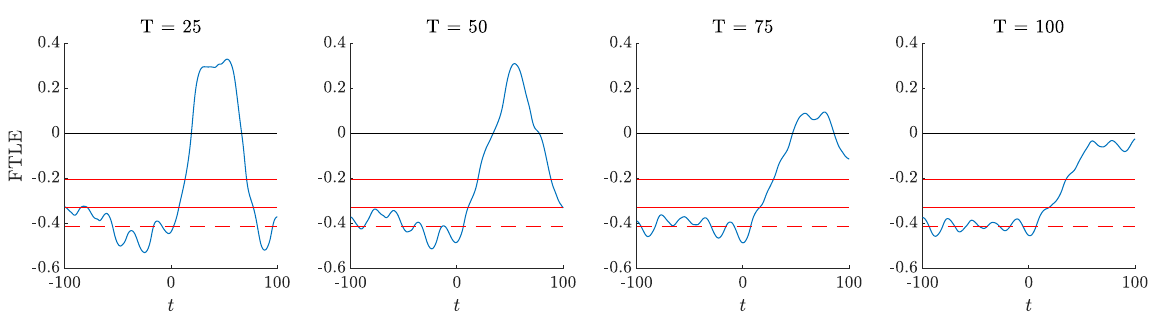}\par
\caption{Approximation of finite-time Lyapunov exponent $t\mapsto\lambda_u(c,T,t)$ of the upper local pullback attractor $u_c$ of \eqref{eq:6Bnumericalproblem}$_c$ for $c=0.999999267212$ and for different values of the length $T$ of the sliding window. In red dashed line, the Lyapunov exponent $L$ of $\tilde u_{\gamma_+}$; and in red solid lines, $\kappa\, L$ for $\kappa\in\{0.5,0.8\}$, two possible warning signals of overgrowth of $\lambda_u(c,T,t)$.}\label{fig:6BFTLEdifferentTs}
\end{figure*}
Assume that the hypotheses \hyperlink{C}{{\rm\textbf{C}}} (resp. \hyperlink{D}{{\rm\textbf{D}}}) hold. As time increases, the long-term behavior of the
locally pullback attractive solution $a_c$ (resp. $l_c$ and $u_c$) is identical 
for all the values of $c$ such that equations \eqref{eq:4Cparametrictransitionequation}$_c$ are in \textsc{Case} \hyperlink{Ac}{{\rm A}} (resp. \textsc{Cases} \hyperlink{Ad}{{\rm A}}, \hyperlink{C1d}{{\rm C1}} or \hyperlink{C2d}{{\rm C2}}):
all of them approach $\tilde a_{\gamma_+}$ (resp. $\tilde l_{\gamma_+}$ and/or $\tilde u_{\gamma_+})$ as time increases. For this reason, the corresponding classical Lyapunov exponents, which only gather information about the asymptotic behaviour of the solutions,
do not provide valuable information
regarding the proximity of $c$ to a tipping point $c_0$ (recall Definition~\ref{def:3.1tippingtracking}).
In this section, we show that, on the contrary,
the finite-time Lyapunov exponents of these solutions are effective (early-warning) signals: as the rate approaches the tipping point $c_0$, there is a moment in time where the FTLE become positive 
for suitable choices of the integration interval.
It is remarkable that this is true despite the fact that equations \eqref{eq:4Cparametrictransitionequation}$_c$ are not necessarily of slow variation.

Let us focus on the concave case, i.e., let us assume the family of hypotheses \hyperlink{C}{{\rm\textbf{C}}}.
Whenever the locally pullback solution $a_c$ is globally defined, we
define its \emph{backward finite-time Lyapunov exponent of length $T>0$ at time $t\in\mathbb{R}$} by
\begin{equation}\label{eq:5Cfinitetimelyapunovexponent}
 \lambda_{a}(c,T,t)=\frac{1}{T}\int^t_{t-T} f_x(s,a_c(s),\Gamma^c(s))\, ds\,.
\end{equation}
We note that, in real-case scenarios, the integrand can be approximated by the rate of separation of close trajectories if a sufficient amount of experimental data for different initial conditions is available.

Let the dynamics of \eqref{eq:4Cparametrictransitionequation}$_{c_0}$ be in \textsc{Case} \hyperlink{Bc}{{\rm B}}, and assume that \eqref{eq:4Cparametrictransitionequation}$_c$
is in \textsc{Case} \hyperlink{Ac}{{\rm A}} for all $c\in(c_1,c_0)$ or all $c\in(c_0,c_2)$.
Let $I$
represent either $(c_1,c_0]$ or $[c_0,c_2)$.
Theorem 4.16 of Ref.~\onlinecite{lno1} can be used to show the continuity of $I\times\mathbb{R}\rightarrow\mathbb{R}$, $(c,t)\mapsto \lambda_a(c,T,t)$ for fixed $T>0$.

Let $t_0\in\mathbb{R}$ be fixed. Since $\lim_{t\rightarrow\infty} (a_{c_0}(t)-\tilde r_{\gamma_+}(t))=0$ and the Lyapunov spectrum of $\tilde r_{\gamma_+}$ (see Subsection~\ref{subsec:2CLyapunovexponents} of Appendix~\ref{appendix1}) is strictly positive, the techniques of Theorems~9.05 and 9.12 of Ref.~\onlinecite{nemytskii1} (which make use of the skewproduct formulation) ensure that
\begin{equation*}
    \lim_{T\rightarrow\infty}\;\inf_{t-T\geq t_0}\frac{1}{T}\int_{t-T}^t f_x(s,a_{c_0}(s),\Gamma^{c_0}(s))\, ds>0\,.
\end{equation*}
Analogously, since $\lim_{t\rightarrow-\infty} (a_{c_0}(t)-\tilde a_{\gamma_-}(t))=0$ and the Lyapunov spectrum of $\tilde a_{\gamma_-}$ is strictly negative,
\begin{equation*}
    \lim_{T\rightarrow\infty}\;\sup_{t\leq t_0}\frac{1}{T}\int_{t-T}^t f_x(s,a_{c_0}(s),\Gamma^{c_0}(s))\, ds<0\,.
\end{equation*}
That is, for our fixed $t_0\in\mathbb{R}$, there exist $T_0>0$ and $\rho>0$ such that, if $T\geq T_0$ and $t\geq t_0+T$, then $\lambda_a(c_0,T,t)\geq\rho$ and, if $T\geq T_0$ and $t\leq t_0$, then $\lambda_a(c_0,T,t)\leq-\rho$.

Therefore, once a suitable $T\geq T_0$ is chosen, the map $t\mapsto \lambda_a(c_0,T,t)$ changes sign from strictly negative values for $t\leq t_0$ to strictly positive values for $t\geq t_0+T$.
Let us now fix
a compact interval $J$ for which $t\mapsto\lambda_a(c_0,T,t)$ is strictly positive, and deduce from the continuity of $I\times J\to\mathbb R,
\;(c,t)\mapsto \lambda_a(c,T,t)$ the existence of $\delta>0$ such that
$\lambda_a(c,T,t)>0$ if $t\in J$ and $c\in I$ satisfies $|c-c_0|<\delta$.
In other words, once chosen a suitable $T\ge T_0$, we get a positive
value of $\lambda_a(c,T,t)$ for a finite interval of time $t\in J$
if $c$ is close enough to $c_0$, despite the fact that \eqref{eq:4Cparametrictransitionequation}$_c$
is in \textsc{Case} \hyperlink{Ac}{{\rm A}} and hence $a_c$ is hyperbolic attractive.
In addition, in general, a larger $T$ forces us to take $c$ closer to $c_0$ in order to have $a_c$ approaching $a_{c_0}$ in the intervals $[t-T,t]$ for $t\in J$. So, the larger $T$ is, the smaller the neighborhood of $c_0$ for which we find positive exponents.

The same arguments show the same result in the d-concave case,
i.e., under the family of hypotheses \hyperlink{D}{{\rm\textbf{D}}}, if \eqref{eq:4Cparametrictransitionequation}$_{c_0}$ is in \textsc{Case} \hyperlink{B1d}{{\rm B1}}
(resp.~\hyperlink{B2d}{{\rm B2}}) and \eqref{eq:4Cparametrictransitionequation}$_{c}$ is either in \textsc{Case} \hyperlink{Ad}{{\rm A}} or in \textsc{Case} \hyperlink{C1d}{{\rm C1}}
(resp. \eqref{eq:4Cparametrictransitionequation}$_{c}$ is either in \textsc{Case} \hyperlink{Ad}{{\rm A}} or in \textsc{Case} \hyperlink{C2d}{{\rm C2}}) for $c\in I$.
Now, we work with the finite-time exponents
$\lambda_l(c,T,t)$ (resp. $\lambda_u(c,T,t)$) defined by replacing
$a_c$ by $l_c$ (resp. by $u_c$) in \eqref{eq:5Cfinitetimelyapunovexponent}.
The proof of this assertion requires to check the properties analogue 
to those of Theorem 4.16 of Ref.~\onlinecite{lno1} for this case, which is straightforward.

\section{Rate-induced tracking phenomena}\label{sec:Vrateinducedtrackingphenomena}
In this section, we present examples of rate-induced tracking.
We say that a parametric family \eqref{eq:4Cparametrictransitionrate}
satisfying hypotheses \hyperlink{C}{{\rm\textbf{C}}} or \hyperlink{D}{{\rm\textbf{D}}} exhibits \emph{rate-induced tracking} if \eqref{eq:4Cparametrictransitionrate}$_c$ is in \textsc{Case} \hyperlink{Ac}{{\rm A}} for sufficiently large values of the rate $c>0$ and in \textsc{Case} \hyperlink{Cc}{{\rm C}} for sufficiently small values of $c>0$.
That is, if the parameter shift takes place at high rate, the system remains in the desired state (\textsc{Case} \hyperlink{Ac}{{\rm A}}), while if it takes place at low rate, a critical transition takes place (\textsc{Case} \hyperlink{Cc}{{\rm C}}).

The concept of rate-induced tracking presented in this section covers a number of models for which a sufficient increase in the transition rate is beneficial to guarantee tracking.

In Subsection~\ref{subsec:5Cconcaverateinducedtracking}, we use $d>0$ for the rate of the transition function $\Delta$, and we also consider the phase of a transition function $\Gamma$ in the sense of \eqref{eq:4Cparametrictransitionphase}.
\begin{figure*}
\centering
\includegraphics[width=\textwidth]{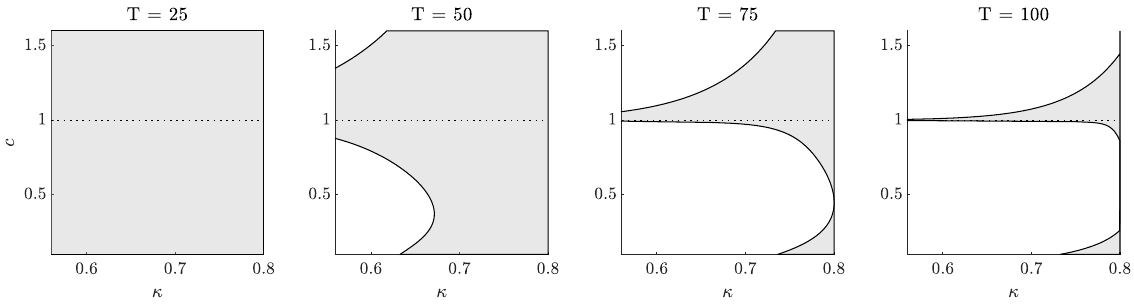}\par
\caption{In grey, for different values of $T>0$, the numerical approximation of the region of pairs $(\kappa,c)$ for which there exists $t\in[-400,400]$ such that $\lambda_u(c,T,t)\geq\kappa\, L$, that is, the grey region contains the points for which the EWS calculated with the finite time-Lyapunov exponents of the locally pullback attractive solution $u_c$ of \eqref{eq:6Bnumericalproblem}$_c$ is issued. In black dashed line, an approximation to the tipping point $c_0$.}\label{fig:6Bdetectionregions}
\end{figure*}

\subsection{Rate-induced tracking in d-concave models}\label{subsec:5Crateinduced}
In this section, we investigate the occurrence of rate-induced tracking in d-concave models subjected to a time-dependent variation of parameter. At first, we focus on the case of variation at constant rate and we numerically show the effectiveness of finite-time Lyapunov exponents to detect the proximity of tipping points.
Afterwards, the case of time-dependent rate is considered, and we present geometric ingredients that help to conclude the occurrence of tipping or tracking at finite time.

\subsubsection{An example with constant rate}\label{subsubsec:VA1anexamplewithtimeindependentrate}
Let us illustrate with an example how rate-induced tracking may arise.
We consider the d-concave model \eqref{eq:3dconcavenoncubicallee} plus a rate parametric migration term $\Gamma(ct)\phi(t)$ with positively bounded-from-below function $\phi$, which describes the arrival of individuals to the population:
\begin{equation}\label{eq:6Bnumericalproblem}
x'=R(t,x)+\Gamma(ct)\phi(t)\,,
\end{equation}
where
\begin{equation}
\label{eq:Rex}    R(t,x)=r(t)\,x\,\left(1-\frac{x}{K(t)}\right)\frac{x-\mu(t)}{\nu(t)+x}\,.
\end{equation}
We recall that $r$, $K$, $\nu$ and $\nu+\mu$ are also positively bounded-from-below. In fact, we also assume that the  functions $r$, $K$, $\mu$, $\nu$ and $\phi$ are quasiperiodic, so that, for any fixed $\gamma\in\mathbb{R}$, the function $f(t,x,\gamma)=R(t,x)+\gamma\,\phi(t)$ is uniformly quasiperiodic for $x$ on compact sets.
Note that the monotonicity of $f$ on the parameter $\gamma$ is essential for the purposes of this section.

In bird populations, some of the factors that produce changes in the size of the migrant population (diseases, storms, adverse winds, orientation errors...) are not persistent over time (especially in species with continental distribution).\cite{rappole1} The natural modeling for these phenomena are parameter shifts with the same asymptotic limits at $+\infty$ and $-\infty$. 
For simplicity, we choose a unimodal function $\Gamma(t)$ to model the amplitude of the migration term: a linear transformation of a Cauchy distribution given by
\begin{equation}\label{eq:pulse}
\Gamma(t)=\gamma_++\frac{\gamma_*-\gamma_+}{1+b\,t^2}\,,
\end{equation}
with $b>0$.
Therefore, $\Gamma^c(t)=\Gamma(ct)$ is our parameter shift. Note that it describes an impulse from the baseline value $\gamma_+$ at $\pm\infty$  to a peak value $\gamma_*$ at $t=0$.
As a consequence, the past and future equations of our system are the same. The parameter $c>0$ represents the speed of the impulse: a greater $c$ corresponds to a smaller window at which $\Gamma(ct)$ is numerically distinguishable from $\gamma_+$.

To get an equation satisfying the family of hypotheses \hyperlink{D}{{\rm\textbf{D}}}, we assume that
\begin{equation}\label{eq:6Bpastandfutureequations}
    x'=R(t,x)+\gamma_+\phi(t)
\end{equation} has three hyperbolic solutions $\tilde l_{\gamma_+}<\tilde m_{\gamma_+}<\tilde u_{\gamma_+}$, which are uniformly separated: as said in Subsection~\ref{subsec:4Bdconcavehypotheses}, the upper and lower ones are hyperbolic attractive while the middle one is hyperbolic repulsive.
Quasiperiodicity guarantess that the system is uniquely ergodic, which, in turn, ensures that the Lyapunov spectra of $\tilde l_{\gamma_+}$, $\tilde m_{\gamma_+}$ and $\tilde u_{\gamma_+}$ reduce to a point each, and hence we will refer to them as to the Lyapunov exponents of $\tilde l_{\gamma_+}$, $\tilde m_{\gamma_+}$ and $\tilde u_{\gamma_+}$ respectively.

In the past and future equations \eqref{eq:6Bpastandfutureequations}, the upper attractive hyperbolic solution $\tilde u_{\gamma_+}$ represents the healthy state of the population, while $\tilde l_{\gamma_+}$, which is near 0, represents extinction.
The upper locally pullback attractive solution $u_c$ of the transition equation \eqref{eq:6Bnumericalproblem}$_c$ represents the evolution of a population, which departs from the healthy state $\tilde u_{\gamma_+}$ in the past. If \eqref{eq:6Bnumericalproblem}$_c$ is in \textsc{Case} \hyperlink{Ad}{{\rm A}}, then $u_c$ approaches $\tilde u_{\gamma_+}$ as time tends to $+\infty$, that is, the population persists, while, if \eqref{eq:6Bnumericalproblem}$_c$ is in \textsc{Case} \hyperlink{C2d}{{\rm C2}}, then $u_c$ approaches $\tilde l_{\gamma_+}$ as time tends to $+\infty$, that is, the population undergoes extinction.
%This is the situation we are interested in, due the biological significance of a critical extinction.

Let us consider the $\gamma$-parametric family of intermediate parameter problems
\begin{equation}\label{eq:6Bstaticproblems}
x'=R(t,x)+\gamma\,\phi(t)\,.
\end{equation}
The arguments proving Theorem~\ref{th:4Dd-concavebifurcationtheorem} can be adapted to check the existence of an open interval $I$ such that 
\eqref{eq:6Bstaticproblems}$_\gamma$ has three hyperbolic solutions for $\gamma\in I$ and just one for $\gamma\not\in I$. 
Reasoning as in Theorem~4.4 of Ref. \onlinecite{dno3}, we check that \eqref{eq:6Bnumericalproblem}$_c$ is in \textsc{Case} \hyperlink{Ad}{{\rm A}} for all $c>0$ if $\Gamma(\mathbb{R})\subset I$. Hence, to investigate the occurrence of tipping, we choose $\Gamma(0)=\gamma_*<\inf I$. 
In this way, $\Gamma(ct)$ (and hence the immigration) decreases from $\gamma_+$ to $\gamma_*$ as $|t|$ decreases.

\begin{figure*}
\centering
\includegraphics[width=\textwidth]{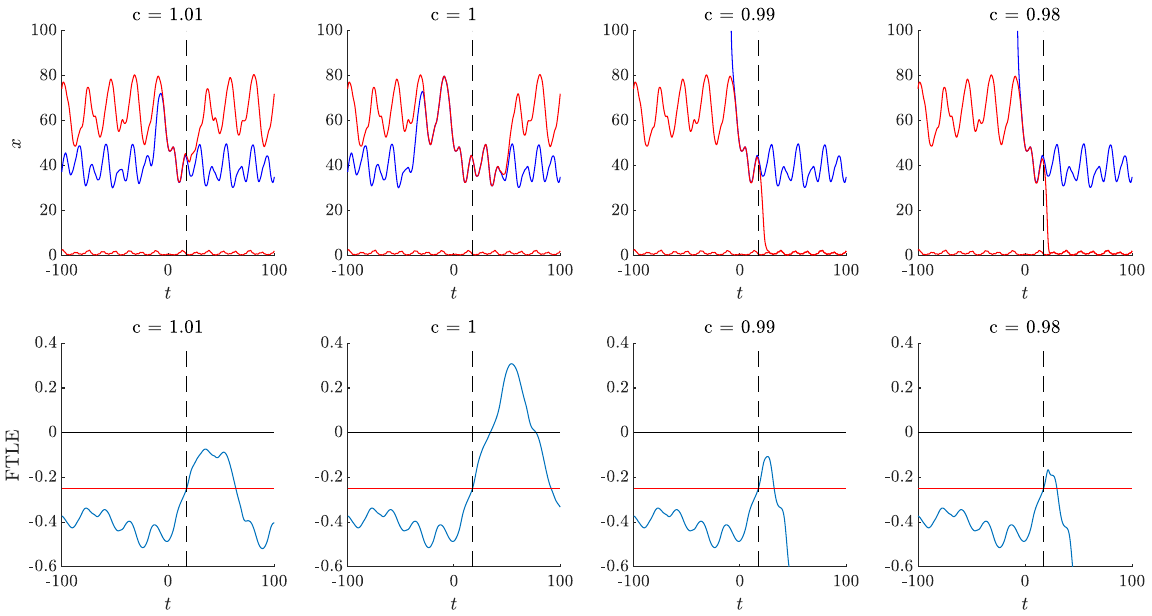}\par
\caption{In the first line of pictures, approximation of the locally pullback attractive solutions $l_c<u_c$ (solid red) and the locally pullback repulsive solution $m_c$ (dashed blue) of \eqref{eq:6Bnumericalproblem}$_c$ for different values of the rate $c>0$.
In the second line, approximation of the finite-time Lyapunov exponent $t\mapsto\lambda_u(c,T,t)$ (FTLE) of $u_c$ for a sliding window of $T=50$ units of time and different values of the rate $c$. The red horizontal line represents the warning threshold $\kappa\,L$ for $\kappa=0.6$. The black dashed vertical line represented in both first and second lines of graphs represents the time at which the warning signal is produced, that is, the first $t\in\mathbb{R}$ such that $\lambda_u(c,T,t)\geq\kappa\, L$.}\label{fig:42simulations}
\end{figure*}
The arguments of Subsection~5.1 of Ref.~\onlinecite{dno3} prove that there exists a unique value of the parameter $c_0>0$ such that \eqref{eq:6Bnumericalproblem}$_c$ is in \textsc{Case} \hyperlink{Ad}{{\rm A}} if $c>c_0$, in \textsc{Case} \hyperlink{B2d}{{\rm B2}} if $c=c_0$, and in \textsc{Case} \hyperlink{C2d}{{\rm C2}} if $c\in(0,c_0)$.
That is, a gradual variation between $\gamma_+$ and $\gamma_*$ means extinction, while a faster one (concentrated around $t=0$) prevents it.

If we consider a phase dependence in \eqref{eq:6Bnumericalproblem} of the form \begin{equation}
    x'=R(t,x)+\Gamma(ct+d)\phi(t)
\end{equation}
with $d\in\mathbb{R}$ for the same $\Gamma$, then we are also in a rate-induced tracking scenario: the same arguments show that it is in \textsc{Case} \hyperlink{Ad}{{\rm A}} for sufficiently large $c>0$ and in \textsc{Case} \hyperlink{C2d}{{\rm C2}} for sufficiently small $c>0$, but the uniqueness of the tipping point is not ensured by the previous arguments.
\par\vspace{2ex}
\textbf{Numerical evidence.} Let us numerically illustrate  the previously described example of rate-induced tracking in a d-concave problem subjected to a time-dependent variation of parameter at a constant rate. Specifically, we aim to showcase the effectiveness of finite-time Lyapunov exponents (FTLE) in detecting the proximity of the unique tipping point $c_0$. To this end, in Figures~\ref{fig:6BFTLEdifferentTs}, \ref{fig:6Bdetectionregions} and \ref{fig:42simulations}, we take
\begin{equation*}
\begin{split}
r(t)=1.5+\sin^2(t/4)&,\quad
K(t)=40+40\sin^2(t\sqrt{5}/16),\\
\mu(t)=30+30\sin^2(t/4)&,\quad
\nu(t)=40+40\sin^2(t\sqrt{5}/16),\\
\phi(t)=0.75+&0.5\sin^2(t\sqrt{5}/2),\\
\gamma_+=1.5,\quad 
\gamma_*=&0.8,\quad 
b=0.02386.
\end{split}
\end{equation*}
In this case, we numerically approximate the Lyapunov exponent $L\approx-0.4134$ of $\tilde u_{\gamma_+}$.

We have explained in Subsection~\ref{subsec:5BFTLE} that, if $T$ is suitably chosen and $c$ is close enough to $c_0$, then the FTLE $t\mapsto \lambda_u(c,T,t)$ of $u_c$ is positive during some periods of time $t$. 
Figure~\ref{fig:6BFTLEdifferentTs} depicts the FTLE calculated at  $c=0.999999267212$ for four different values of integration length $T$: there is numerical evidence of \textsc{Case} \hyperlink{Ad}{{\rm A}} for \eqref{eq:6Bnumericalproblem}$_c$, and that this is no longer the case for $c=0.999999267211$: hence, the tipping point $c_0$ must lie very close to these values. Let us briefly discuss the behaviour of the FTLE depending on the length of the interval of integration.
On the one hand, we note that the four lengths of integration intervals $T=25, 50, 75, 100$, are such that the FTLE $\lambda_u(c,T,t)$ approximate the Lyapunov exponent $L$ reasonably well whenever $-t$ is large enough. On the other hand, any of the lengths $T=25,50,75$ is small enough to ensure that the FTLE exceed zero at some point in time.
Note that values of $T$ smaller than $25$ substantially worsen the approximation of the Lyapunov exponents of past and future, while a positive FTLE for $T=100$, is obtained only if the rate approximates $c_0$ beyond the eleventh decimal digit.

Since the range of values of $c$ for which we find positive FTLE may be too narrow to react in time and prevent tipping, we define more sensitive indicators of the proximity of a tipping point: instead of waiting for the change of sign in the FTLE, we focus on the overgrowth of the values of FTLE anticipating such change of sign.
Recall that for the considered model \eqref{eq:6Bnumericalproblem}, with $R(t,x)$ and $\Gamma(t)$ respectively as in \eqref{eq:Rex} and \eqref{eq:pulse}, the Lyapunov exponents $L$ of $\tilde u_{\gamma_+}$ and those of $\tilde u_{\gamma_-}$ are the same (in fact, $\gamma_-=\gamma_+$). Therefore, we can define a suitable alert threshold as $\kappa\, L$, i.e.~a certain fraction $\kappa\in[0,1)$ of the value of the Lyapunov exponent $L$ of $\tilde u_{\gamma_+}$.
Then, an early-warning signal (EWS) is issued the first time $t\in\mathbb{R}$ that $\lambda_u(c,T,t)
\geq \kappa\, L$.
In a more general scenario where $\tilde u_{\gamma_-}$ and $\tilde u_{\gamma_+}$ have different Lyapunov exponents, a definition of EWS taking into account both of them would be beneficial.

In Figure~\ref{fig:6Bdetectionregions}, we find the regions of pairs $(\kappa,c)$ for which there exist $t\in\mathbb{R}$ such that $\lambda_u(c,T,t)\geq \kappa\,L$ for different values of the integration length $T>0$ of the FTLE.
The reduction of the sensitivity of FTLE, which is represented by the decrease of the curves $t\mapsto \lambda_u(0.999999267212,T,t)$ as $T$ increases in Figure~\ref{fig:6BFTLEdifferentTs} is again accurately represented in Figure~\ref{fig:6Bdetectionregions}: the bigger $T$ is, the smaller the shadowed region becomes.
In any of the graphs of Figure~\ref{fig:6Bdetectionregions}, the vertical section of the grey regions over each fixed value of $\kappa$ represents the interval of values of the rate $c>0$ for which EWS are issued.
To avoid the over-detection which apparently presents the case of $T=25$ and the possibly too narrow detection region for rates $c<c_0$ of the cases $T\in\{75,100\}$, we will take $T=50$ for the following numerical experiments.

Figure \ref{fig:6Bdetectionregions} shows another interesting phenomenon: EWS are also detected for $c$ close to 0, which is far away from the tipping point $c_0$. This means that, in this case, the FTLE approximate $0$ for sufficiently small rate. In turn, this fact means slower dynamics for some intervals of time.
The numerical results show that, in addition, these intervals of slow dynamics correspond to values of $t$ for which $\Gamma(ct)$ is close to $\gamma_1$ (which is the first bifurcation value of \eqref{eq:6Bstaticproblems}$_{\gamma}$). (This is not depicted in Figure \ref{fig:6Bdetectionregions}).
The underlying phenomenon is that, for these intervals of time, the dynamics of \eqref{eq:6Bnumericalproblem}$_c$ is quite similar to that of \eqref{eq:6Bstaticproblems}$_{\gamma_1}$. %, which typically has upper Lyapunov exponent equal to 0. 
Clearly, these intervals are longer when $c>0$ decreases.

To some extent, this type of warning (the one for small $c>0$) corresponds to the EWS character of FTLE which have been already considered in Refs. \onlinecite{lenton1,das,moore,schefferbas,dakos,jager1}.
For each $\gamma\in\Gamma(\mathbb R)$, a slow parameter shift causes the dynamics of \eqref{eq:6Bnumericalproblem}$_{c}$ to be very similar to that of \eqref{eq:6Bstaticproblems}$_{\gamma}$ for a ``long time": the time during which $\Gamma(ct)$ is very close to $\gamma$.
For large negative values of the time, the value of $\gamma$ approximated by $\Gamma(ct)$ is close to $\gamma_+$. As time increases, the value of $\gamma$ also decreases towards $\gamma_*$. As this value of $\gamma$ approaches $\gamma_1$ from the right, the short-term dynamics of \eqref{eq:6Bnumericalproblem}$_c$ gets slower, and this slowing down is reflected by FTLE closer to 0. And the smaller $c>0$ is, the more pronounced the slowing down is.

To sum up, Figure \ref{fig:6Bdetectionregions} shows two phenomena producing EWS of different nature. On the one hand, the detection of proximity of \textsc{Case} \hyperlink{Bd}{{\rm B}} on a neighborhood of $c_0$, and, on the other hand, a critical slowing-down for sufficiently small rates $c$ which indicates the moment at which the parameter shift $\Gamma$ escapes from the interval $I$.

In Figure~\ref{fig:42simulations}, the locally pullback attractive solutions $l_c<u_c$ and the locally pullback repulsive solution $m_c$ of \eqref{eq:6Bnumericalproblem}$_c$ are represented together with the FTLE $t\mapsto\lambda_u(c,T,t)$ for $T=50$ and for different values of the rate $c$.
The threshold of warning signals is fixed for the value $\kappa=0.6$.
Recall that the tipping point $c_0$ lies between $c=0.99$ and $c=1.00$.
As it can also be seen in Figure~\ref{fig:6Bdetectionregions}, in this example, the tipping point is better detected from the region $c>c_0$ of \textsc{Case} \hyperlink{Ad}{{\rm A}}: $c=0.99$ is closer to the tipping point $c_0$ than $c=1.01$ and, however, the overgrowth of FTLE is better observed for $c=1.01$.
The locally pullback attractive solution $u_c$ and the locally pullback repulsive solution $m_c$ approach each other in compact intervals of time as $c\rightarrow c_0$.
The vertical dashed line represents the value of time at which the warning signal is triggered.
In general, the closer we are to $c_0$, the earlier in time we obtain a warning signal.
For $c=0.99$ and $c=0.98$, we notice that the warning signal is issued an interval of time before the population gets extincted.
Although this interval of time can seem small in the presented graphs, we recall that the scale of time at which the transition takes place (maybe celestial or climatic time scales) can be 
larger than the time needed to reverse tipping by some kind of external action, as we will see in the numerical examples of Subsection~\ref{subsec:6Cnoreturnrateinducetracking}: the interval of time which goes from the warning signal to the extinction of the species is the one on which we should engage actions to prevent tipping.
\subsubsection{Geometric ingredients with time-dependent rate}\label{subsubsec:VA2geometricingredients}
\begin{figure}
\centering
\includegraphics[width=0.45\textwidth]{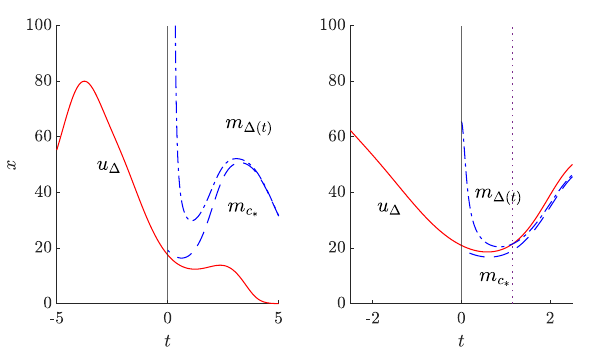}\par
\caption{Illustration of no-return points $u_\Delta(s_3)<m_{c_*}(s_3)$, for any $s_3\geq0$, (left) and safe points $u_\Delta(s_2)>m_{\Delta(s_2)}(s_2)$, for any $s_2>1.15$ (right). In red solid line, the locally pullback attractive solution $u_\Delta$ of \eqref{eq:6Cnonconstant}, in blue dashed-dotted line, the curve of values of the static \eqref{eq:6constantproblems}$_{\Delta(t)}$ repellers $t\mapsto m_{\Delta(t)}(t)$ (upper), and, in blue dashed line, the locally pullback repulsive solution of the future $m_{c_*}$ (lower). For both pictures, $f(t,x,\gamma)=r(t)\,x\,(1-x/K(t))-(52-13\gamma)\,x\,/(x+10)$ with $r(t)=2+\sin(t)$ and $K(t)=90+18\sin^2(t\sqrt{5}/2)$, and $\Gamma(t)=-550/(1000+t^2)$. In both cases, $c=19.652326446580$ is the nearest point of \textsc{Case} {\rm\textbf{A}} with 12 decimal digits to $c_0$, $t_0=0$ and the warning point $s_1\approx-24.34$ does not appear in the representation. On the left, $\Delta(t)=35-300/(10+t^2)$, on the right $\Delta(t)=30-147/(6+4t^2)$. (Remark 4.7 of Ref.~\onlinecite{dno3} shows that the strictly increasing character of $\gamma\mapsto f(t,x,\gamma)$ for $(t,x)\in\mathbb{R}\times(0,\infty)$ is enough).}\label{fig:6Cillustration}
\end{figure}
Under some assumptions on the monotonicity of the parametric variation of the right-hand side function, we explain some geometric ingredients of rate-induced tracking for a problem with time-dependent rate.
We show the existence of a curve $t\mapsto m_{\Delta(t)}(t)$ (resp. $t\mapsto m_{c_*}(t)$) such that if the upper locally pullback attractive solution is above (resp. below) that curve at some time, then tracking (resp. tipping) takes place.
%There exist curves that, once they are traversed by the locally pullback attractive solution, tell us whether the transition is of tipping or tracking.

Let $\Gamma\colon\mathbb{R}\rightarrow\mathbb{R}$ and $f\colon\mathbb{R}\times\mathbb{R}\times\mathbb{R}\rightarrow\mathbb{R}$ be such that the pair $(f,\Gamma)$ satisfies the family of d-concave hypotheses \hyperlink{D}{{\rm\textbf{D}}}, that $\gamma\mapsto f(t,x,\gamma)$ is strictly increasing for all $(t,x)\in\mathbb{R}\times\mathbb{R}$, and that $\Gamma$ is nonincreasing on $(-\infty,0]$ and nondecreasing on $[0,\infty)$, with $\gamma_-=\gamma_+$.
The parameter shift $\Gamma$ represents an impulse which takes its maximum value at $t=0$.
Notice that, in this case, $c\mapsto f(t,x,\Gamma(ct))$ is nondecreasing for all $(t,x)\in\mathbb{R}\times\mathbb{R}$.
Assume also that we are in a rate-induced tracking problem for the upper locally pullback attractive solution with a unique tipping point, that is, there exists $c_0>0$ such that 
\begin{equation}\label{eq:6constantproblems}
x'=f(t,x,\Gamma(ct))
\end{equation}
is in \textsc{Case} \hyperlink{Ad}{{\rm A}} if $c>c_0$, in \textsc{Case} \hyperlink{B2d}{{\rm B2}} if $c=c_0$, and in \textsc{Case} \hyperlink{C2d}{{\rm C2}} if $c\in(0,c_0)$. (Theorem 4.8 of Ref.~\onlinecite{dno3} shows that this is a natural situation under the previous assumptions).
As previously, for $c>0$, let $l_c$ and $u_c$ be the locally pullback attractive solutions of \eqref{eq:6constantproblems}$_c$, let $m_c$ the locally pullback repulsive solution of \eqref{eq:6constantproblems}$_c$, and let $t\mapsto x_c(t,s,x_0)$ denote the maximal solution of \eqref{eq:6constantproblems}$_c$ with initial condition $x_0$ at time $s\in\R$.
%Notice that, under the assumed hypotheses, the nondecreasing character of $c\mapsto f(t,x,\Gamma(ct))$ ensures that $c\mapsto l_c(t),u_c(t)$ are nondecreasing for all $t\in\mathbb{R}$ and that $c\mapsto m_c(t)$ is nonincreasing for all $t\in\mathbb{R}$.\cite{dno1}

Consider a continuous function $\Delta\colon\mathbb{R}\rightarrow(0,\infty)$  such that $\lim_{t\rightarrow\pm\infty} \Delta(t)=c_*$ for some $c_*>0$ and $c_0<c_*$. Furthermore, assume that there exists $t_0\in\mathbb{R}$ such that $\Delta$ is nondecreasing on $[t_0,\infty)$.
We also assume $|\Delta'|$ to be sufficiently small.

We study the transition equation with time-dependent rate
\begin{equation}\label{eq:6Cnonconstant}
x'=f\big(t,x,\Gamma\big(\Delta(t)\,t\big)\big)\,.
\end{equation}
Understanding $\Delta$ as a new parameter shift, we can apply the theory presented in Subsection \ref{subsec:4Bdconcavehypotheses} to \eqref{eq:6Cnonconstant}.
Therefore, there exist locally pullback attractive solutions $l_\Delta$ and $u_\Delta$ and a locally pullback repulsive solution $m_\Delta$ of \eqref{eq:6Cnonconstant} which govern the dynamics of the equation.
We will denote the maximal solutions of \eqref{eq:6Cnonconstant} by $t\mapsto x_\Delta(t,s,x_0)$.

Now, we will describe the geometric ingredients which will help us determine in which dynamical case \eqref{eq:6Cnonconstant} lies.
First, let $c_1=\inf\Delta$, if $c_0<c_1$, then \eqref{eq:6constantproblems}$_{c_1}$ is in \textsc{Case} \hyperlink{Ad}{{\rm A}}: since $f(t,x,\Gamma(\Delta(t)\,t))\geq f(t,x,\Gamma(c_1\,t))$ for all $t\in\mathbb{R}$, the order relation on the right-hand sides ensures the order relation on the bounds of the set of bounded solutions,\cite{dno1} so $u_\Delta(t)\geq u_{c_1}(t)$ for all $t\in\mathbb{R}$, and therefore \eqref{eq:6Cnonconstant} is in \textsc{Case} \hyperlink{Ad}{{\rm A}} (see Lemma 4.5(i) of Ref.~\onlinecite{dno1}).
That is, if the parameter shift $\Delta$ does not reach the tipping point $c_0$ of \eqref{eq:6constantproblems}, then there exists no possibility of tipping in \eqref{eq:6Cnonconstant}.
In view of the above, in the case of $\inf\Delta<c_0$, we define the \emph{warning point} $s_1$ as the first time that $\Delta (s_1)=c_0$. It indicates that we are traversing the threshold of the ``safe zone'' $(c_0,\infty)$.
%This result can also be seen as a size-induced tracking problem: once fixed $\Delta$, the equation $x'=f(t,x,\Gamma((c_*+d\,(\Delta(ct)-c_*))\,t)$ is in \textsc{Case} \hyperlink{Ad}{{\rm\textbf{A}}} for sufficiently small values of $d>0$.

Once the parametric variation $\Delta$ takes our system \eqref{eq:6Cnonconstant} out of the ``safe zone'', we shall say that $s_2\geq t_0$ is a \emph{safe point} if $u_\Delta(s_2)>m_{\Delta(s_2)}(s_2)$.
In this case, for any $x_0\in(m_{\Delta(s_2)}(s_2),u_\Delta(s_2))$ it is satisfied that $\lim_{t\rightarrow\infty}(x_{\Delta(s_2)}(t,s_2,x_0)-u_{\gamma_+}(t))=0$, and since a standard comparison argument (using that $f(t,x,\Gamma(\Delta(s_2)\,t))\leq f(t,x,\Gamma(\Delta(t)\,t))$ for $t\geq s_2$) states that $u_\Delta(t)>x_\Delta(t,s_2,x_0)\geq x_{\Delta(s_2)}(t,s_2,x_0)$ for all $t\geq s_2$, then \eqref{eq:6Cnonconstant} is in \textsc{Case} \hyperlink{Ad}{{\rm A}} (see again Lemma 4.5(i) of Ref.~\onlinecite{dno1}).
That is, once we are in the time interval $[t_0,\infty)$ on which $\Delta$ is nondecreasing, if the locally pullback attractive solution of \eqref{eq:6Cnonconstant} reaches a level above the locally pullback repulsive solution of the underlying fixed $c$-parameter problem at this time, then tracking is guaranteed.
Notice that the curve $t\mapsto m_{\Delta(t)}(t)$ is not a solution of any of the previously specified equations, since it takes values of the locally pullback repulsive solutions of different problems of \eqref{eq:6constantproblems} at each time.
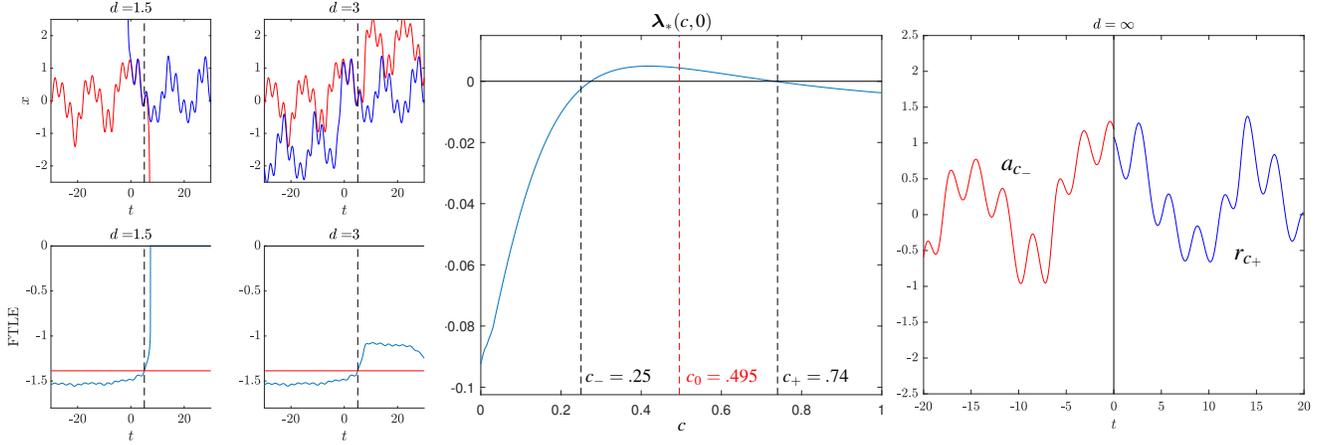
\begin{figure*}
\centering
\includegraphics[trim={0.5cm 0cm 0.7cm 0cm},clip,width=0.32\textwidth]{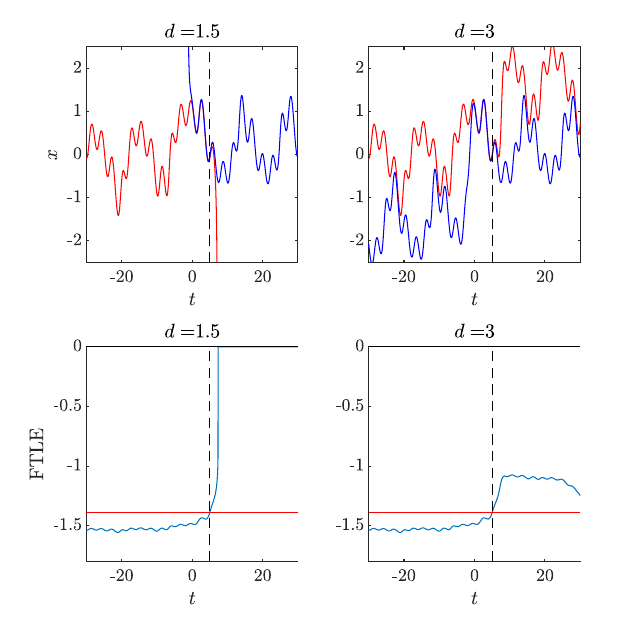}
\begin{overpic}[trim={0.8cm -0.45cm 1.2cm 0.8cm},clip,width=0.325\textwidth]{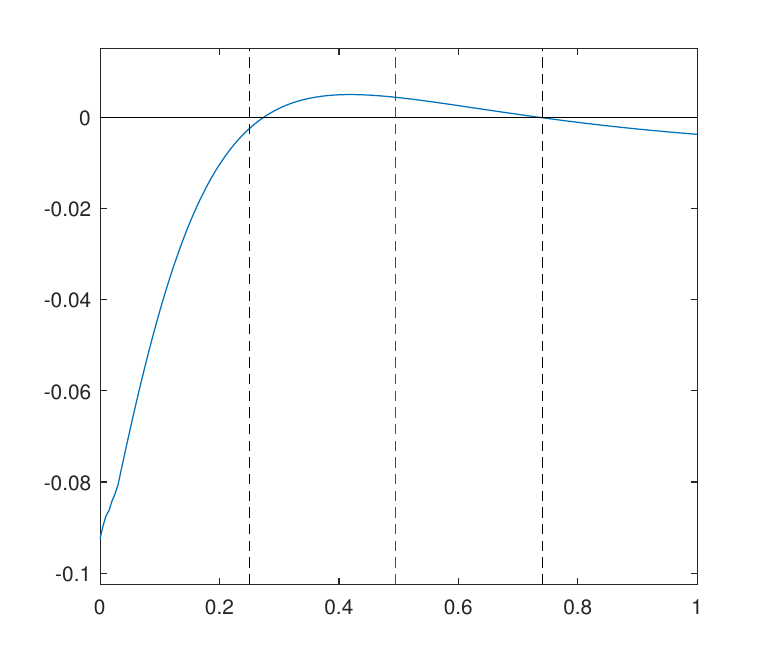}
\put(53,7){\scriptsize $c$}
\put(47,99){\scriptsize $\boldsymbol\lambda_*(c,0)$}
\put(32,18){\scriptsize $c_-=.25$}
\put(77,18){\scriptsize $c_+=.74$}
\put(55,18){\scriptsize \color{red} $c_0=.495$}
\end{overpic}
\begin{overpic}[trim={0.8cm -0.5cm 0.8cm 0.3cm},clip,width=0.32\textwidth]{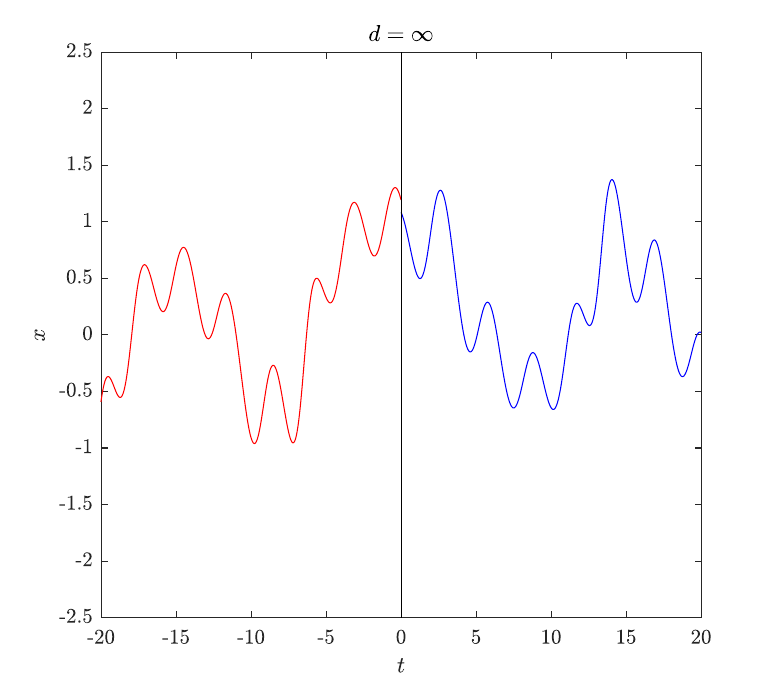}
\put(25,65){$a_{c_-}$}
\put(77,45){$r_{c_+}$}
\end{overpic}
\caption{Rate-induced tracking with time-dependent rate for the concave logistic model \eqref{eq:concave_logistic} with $\Gamma\big(\Delta(dt)\,t\big)=2/\pi \arctan\big(\Delta(dt)\,t\big)$ and $I(t)=-\sin(t/2)-\sin(\sqrt5t)+0.895$, and $\Delta(t)=c_-/(1+e^t)+c_+/(1+e^{-t})$, where $c_-=0.25$ and $c_+=0.74$. On the left-hand side there are four panels. The upper two panels show the locally pullback attracting solution $a_\Delta(t)$ in red, and the locally pullback repelling solution $r_\Delta(t)$ in blue in the extended phase space for the values of $d$ equal to $d=1.5$ (left) and $d=3$ (right). 
The lower panels show the finite-time Lyapunov exponent $t\mapsto \lambda_a(c,T,t)$ (in blue) calculated along the attractor $a_\Delta(t)$ for a sliding window of $T=25$ units of time. In analogy with Figure \ref{fig:42simulations} the red horizontal line represents the warning threshold corresponding to $\kappa=0.9$ of the Lyapunov exponent of the attractor in the past-limit problem. 
In all the first four panels, the black dashed vertical line represents the time at which the warning signal is produced. 
The central panel shows a numerical approximation of the bifurcation map $\boldsymbol\lambda_*(c,0)$ of \eqref{eq:concave_logistic} for $c\in[0,1]$; notice that between the chosen $c_-=0.25$ and $c_+=0.74$ there is a whole interval of values for which $\boldsymbol\lambda_*(c)>0$. The right-hand side panel shows the relative position between the attractor $a_{c_-}$ and the repeller $r_{c_+}$ at $t=0$ providing the condition of tracking for sufficiently high rate $d$ thanks to analysis of the switching problem at $t=0$.
The change of variables $y=x+\gamma$, for $\gamma>\!\!>0$, transforms the coefficients and hyperbolic solutions into others with biological meaning.}
\label{fig:variable_rates}
\end{figure*}

Also outside the ``safe zone'', we shall say that $s_3\geq t_0$ is a \emph{no-return point} if $u_\Delta(s_3)<m_{c_*}(s_3)$.
An analogous argument to the previous one shows that then \eqref{eq:6Cnonconstant} is in \textsc{Case} \hyperlink{C2d}{{\rm C2}}.
That is, once the nondecreasing character of $\Delta$ has started, if the locally pullback attractive solution of \eqref{eq:6Cnonconstant} reaches a level below the locally pullback repulsive solution of the future (and past) equation \eqref{eq:6constantproblems}$_{c_*}$, then the system is irretrievably led to tipping. 
In Figure \ref{fig:6Cillustration}, we can find two numerical examples of these points. In the one on the left, $u_\Delta(0)<m_{c^*}(0)$ indicates that tipping is going to occur, while in the one on the right $u_\Delta(s_2)>m_{\Delta(s_2)}(s_2)$ for $s_2>1.15$ points out tracking.

These safe and no-return points contribute to clarify the understanding of tipping in d-concave systems as finite-time phenomena: there is no need to calculate the dynamics of the locally pullback attractive solutions along the whole real line to determine if we are in a situation of tipping or tracking. This idea was already clear in the concave case: the blow-up of the solutions in \textsc{Case} \hyperlink{Cc}{{\rm C}} takes place at finite time.\cite{lno3,kuehnlongo}

Symmetric definitions could be done in the case we are interested in the dynamics of $l_\Delta$.
Hereby, we have analyzed the case of $\Delta$ nondecreasing.
The same treatment can be carried out for other cases.
\subsection{Rate-induced tracking in concave models for time-dependent rates and phases}\label{subsec:5Cconcaverateinducedtracking}

\begin{figure*}
\centering
\centering
\includegraphics[trim={0.5cm 0cm 0.7cm 0cm},clip,width=0.32\textwidth]{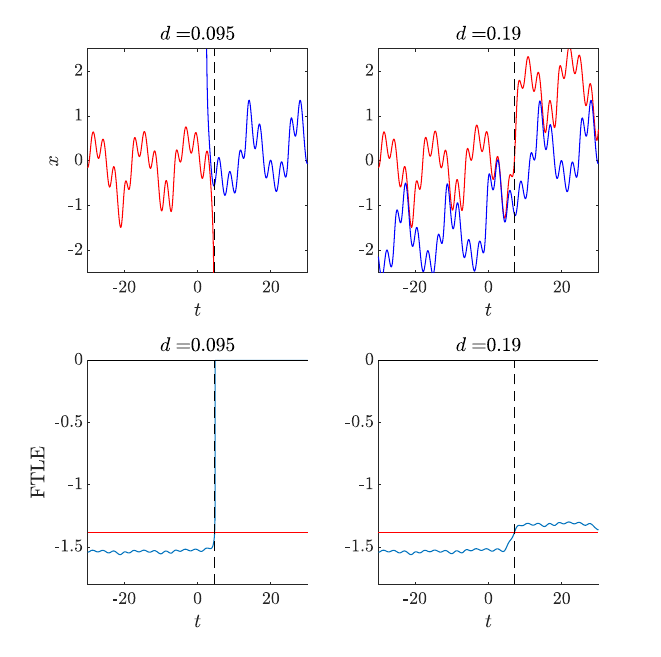}
\begin{overpic}[trim={0.8cm -0.5cm 1.2cm 0.8cm},clip,width=0.312\textwidth]{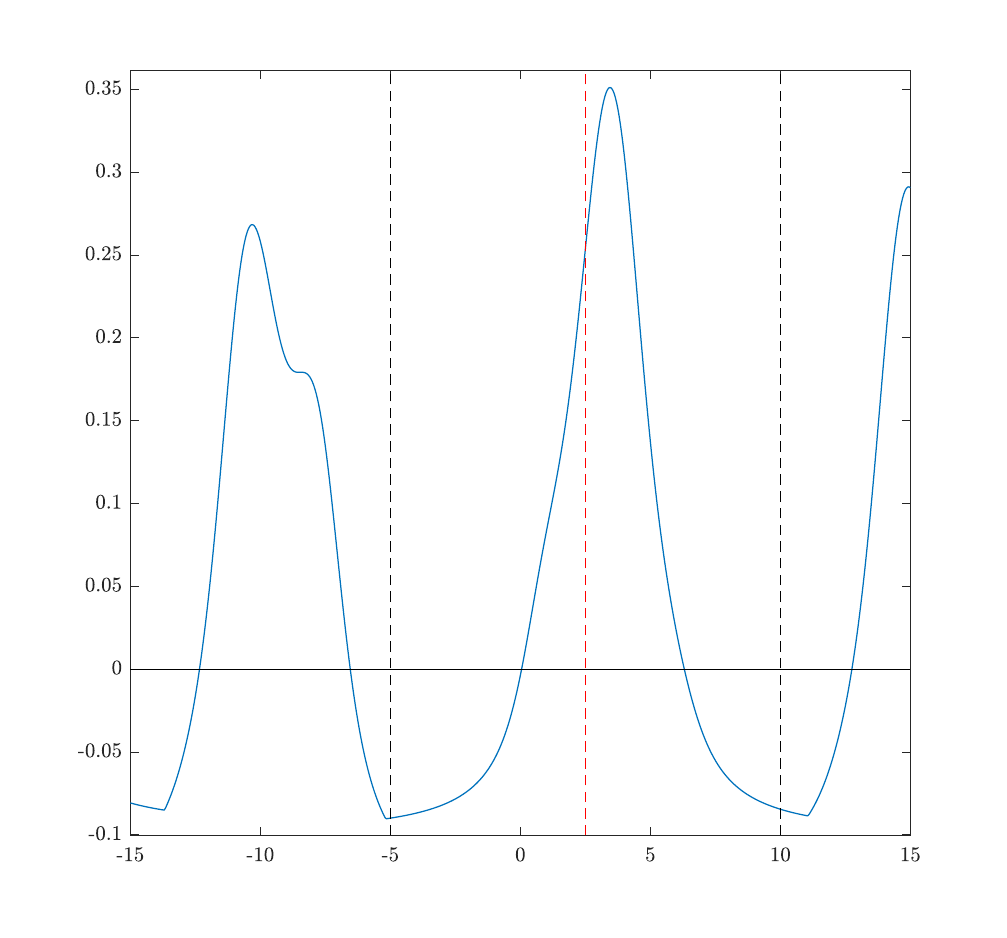}
\put(50,7){\scriptsize $s$}
\put(45,100){\scriptsize $\boldsymbol\lambda_*(1,s)$}
\put(19,92){\scriptsize $s_-=-5$}
\put(63,92){\scriptsize $s_+=10$}
\put(41,92){\scriptsize \color{red} $s_0=2.5$}
\end{overpic}
\begin{overpic}[trim={0.8cm -0.5cm 0.8cm 0.3cm},clip,width=0.32\textwidth]{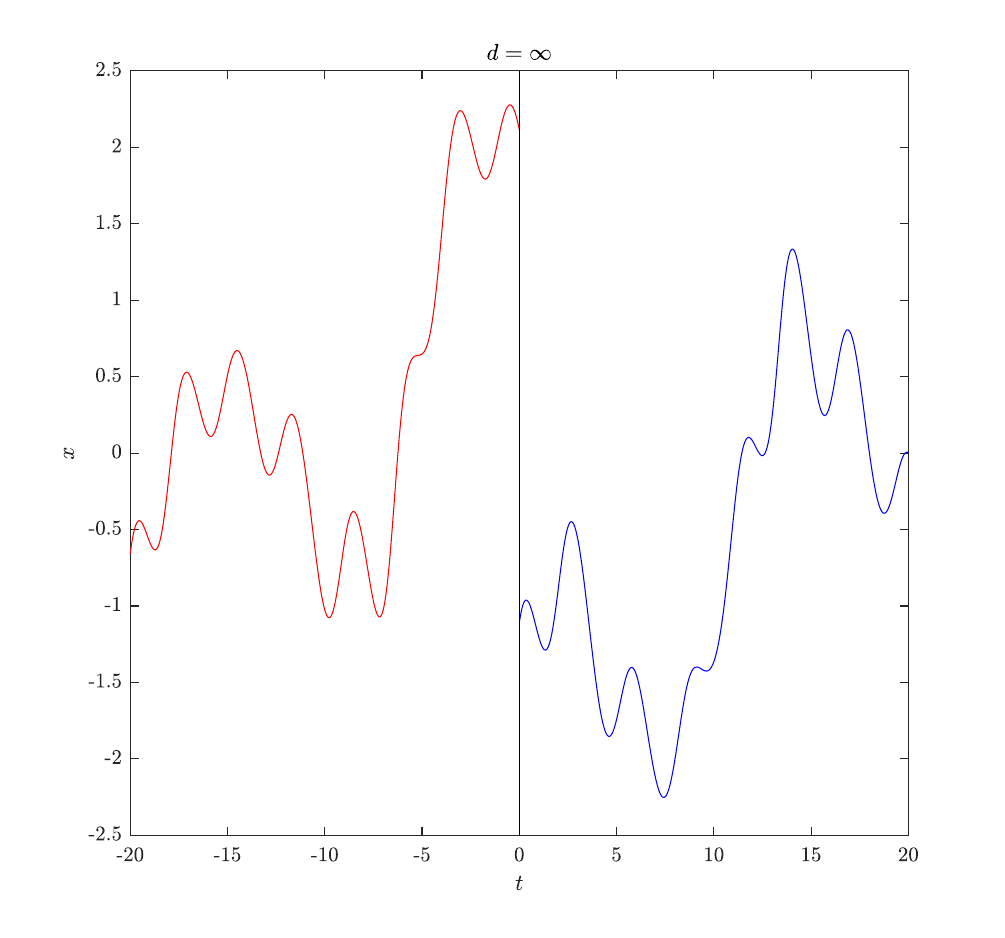}
\put(25,65){$a_{s_-}$}
\put(77,45){$r_{s_+}$}
\end{overpic}
\caption{Rate-induced tracking with time-dependent phase for the concave logistic model \eqref{eq:concave_logistic} with $\Gamma\big(c(t+\Delta(dt)\big)=(2/\pi)\arctan\big(c(t+\Delta(dt)\big)$, with $c=1$, and $I(t)=-\sin(t/2)-\sin(\sqrt5t)+0.895$, and $\Delta(t)=s_-/(1+e^t)+s_+/(1+e^{-t})$, where $s_-=-5$ and $s_+=10$. On the left-hand side there are four panels. The upper two panels show the locally pullback attracting solution $a_\Delta(t)$ in red, and the locally pullback repelling solution $r_\Delta(t)$ in blue in the extended phase space for the values of $d$ equal to $d=0.095$ (left) and $d=0.19$ (right). 
The lower panels show the finite-time Lyapunov exponent $t\mapsto\lambda_a(c,T,t)$ (in blue) calculated along the attractor $a_\Delta(t)$ for a sliding window of $T=25$ units of time. In analogy with Figure \ref{fig:42simulations} the red horizontal line represents the warning threshold corresponding to $\kappa=0.9$ of the Lyapunov exponent of the attractor in the past-limit problem. 
In all the first four panels, the black dashed vertical line represents the time at which the warning signal is produced. 
The central panel shows a numerical approximation of the bifurcation map $\boldsymbol\lambda_*(10,s)$ of \eqref{eq:concave_logistic} allowing to appreciate that between the chosen $s_-=-5$ and $s_+=10$ there is a whole interval of values for which $\boldsymbol\lambda_*$ is positive. The right-hand side panel shows the relative position between the attractor $a_{s_-}$ and the repeller $r_{s_+}$ at $t=0$ providing the condition of tracking for sufficiently high rate $d$ thanks to analysis of the switching problem at $t=0$.
The change of variables $y=x+\gamma$, for $\gamma>\!\!>0$, transforms the coefficients and hyperbolic solutions into others with biological meaning.}
\label{fig:variable_phases}
\end{figure*}
In this section we aim to show alternative mechanisms producing rate-induced tracking for the simplest concave biological model, the logistic equation \eqref{eq:3concavequadratic} with a migration term.
Specifically, we consider a parameter shift $\widehat\Gamma$ and
\begin{equation}\label{eq:concave_logistic}
x'=-r(t)\,\big(x-\widehat\Gamma(t)\big)^2+I(t)\,,
\end{equation}
which can be rewritten in the form of \eqref{eq:3concavequadratic} as 
\begin{equation}
    x'=-r(t)\, x\,(x-2\widehat\Gamma(t))+\widehat I(t)\,,
\end{equation}
where we have renamed the migration term $\widehat I(t):=I(t)-r(t)\,\widehat\Gamma^2(t)$.
For the sake of simplicity, we shall assume $r(t)\equiv 1$ from now on.
\begin{comment}
Specifically, we shall assume that, for some fixed units, the magnitude of the reproduction rate equals that of the carrying capacity at each instant of time, that is, $r(t)=K(t)=2\widehat\Gamma(t)$ for all $t\in\R$. 
Then, \eqref{eq:3concavequadratic} reads as 
\begin{equation}\label{eq:concave_logistic}
\begin{split}
x'&=-2\widehat\Gamma(t)x\big(1-x/\widehat\Gamma(t)\big)\!+ \widehat I(t)%\\
%&=-2\Gamma(t)x\big(1-x/\Gamma(t)\big)-\Gamma^2(t)+ I(t)
=-\big(x-\widehat\Gamma(t)\big)^2\!\!\!+I(t), 
\end{split}
\end{equation}
where we have renamed $ I(t):=\widehat I(t)+\widehat\Gamma^2(t)$ and $\widehat\Gamma(t)$ will be our parameter shift. 
\end{comment}
Note that if the family of hypotheses \hyperlink{C}{{\rm C}} holds for this problem, then 
\begin{equation}\label{eq:5Batractorrepeller}
    x'=-(x-\gamma)^2+I(t)
\end{equation}
has an attractor-repeller pair for every $\gamma\in\R$: the change of variable $y=x-\gamma$ shows that the attractor-repeller pair is given by $\widetilde a(t)-\gamma$ and $\widetilde r(t)-\gamma$, where $\widetilde a(t)$ and $\widetilde r(t)$ are respectively the attractor and the repeller of \eqref{eq:5Batractorrepeller}$_0$. 
Hence, the argument presented in Subection \ref{subsubsec:VA1anexamplewithtimeindependentrate} is not applicable to generate a rate-induced tracking scenario, for there is no interval of values $\gamma$ where the equation \eqref{eq:5Batractorrepeller} is in \textsc{Case} \hyperlink{C}{{\rm C}}.
Therefore, we focus on the case of time-dependent rate and of time-dependent phase. 
Before doing so, let us recall that the results presented in Refs.~\onlinecite{lno1, lno2} allow to construct a bounded and continuous bifurcation map $\boldsymbol\lambda_*(c,s)$ for \eqref{eq:concave_logistic} when a rate $c>0$ and a phase $s\in\R$ intervene in the transition, i.e.~$\widehat\Gamma(t)=\Gamma(c(t+s))$. The sign of the bifurcation map characterizes the dynamical scenario of the transition equation: \textsc{Case} \hyperlink{Ac}{{\rm A}} occurs if and only if $\boldsymbol\lambda_*(c,s)<0$;
\textsc{Case} \hyperlink{Bc}{{\rm B}} occurs if and only if $\boldsymbol\lambda_*(c,s)=0$;
\textsc{Case} \hyperlink{Cc}{{\rm C}} occurs if and only if $\boldsymbol\lambda_*(c,s)>0$. An example of bifurcation map for a problem like \eqref{eq:concave_logistic} can be seen in Figure \ref{fig:variable_rates}.

\begin{comment}

\textsc{Case} \hyperlink{Ac}{{\rm A}} appears if and only if $\widetilde a(t_0)-\gamma_->\widetilde r(t_0)-\gamma_+$;
\textsc{Case} \hyperlink{Bc}{{\rm B}} appears if and only if $\widetilde a(t_0)-\gamma_-=\widetilde r(t_0)-\gamma_+$;
\textsc{Case} \hyperlink{Cc}{{\rm C}} appears if and only if $\widetilde a(t_0)-\gamma_-<\widetilde r(t_0)-\gamma_+$. Most importantly, Theorem 4.5 in Ref.~\onlinecite{lno2} guarantees that if the switching problem is either in \textsc{Case} \hyperlink{Ac}{{\rm A}} or \textsc{Case} \hyperlink{Cc}{{\rm C}} then there is a minimum $c_M>0$ such that \eqref{eq:concave_logistic} showcases the same dynamics for every $c> c_M$. Since we aim to construct a scenario of rate-induced tracking, we will choose $\gamma_-$ and $\gamma_+$ such that $\widetilde a(t)-\gamma_->\widetilde r(t)-\gamma_+$ for all $t\in\R$, i.e.~such that the switching problem at any $t_0\in\R$ is in \textsc{Case} \hyperlink{Ac}{{\rm A}}.

\end{comment}

\par\smallskip
\textbf{Time-dependent rate.}
As for Section \ref{subsubsec:VA2geometricingredients}, we shall consider the time-dependent rate version of \eqref{eq:concave_logistic}, i.e.,
\begin{equation}\label{eq:concave-variable-rate}
x'=-\big[x-\Gamma\big(\Delta(d\,t)\,t\big)\big]^2+I(t),
\end{equation}
with $\Gamma:\R\to\R$ and $\Delta: \R\to(0,\infty)$ continuous functions respectively limiting to $\gamma_-$ and $c_-$ as $t\to-\infty$ and to $\gamma_+$ and $c_+$ as $t\to\infty$. The function $\Delta(t)$ produces the time-dependent rate and the new parameter $d>0$ is the "rate" for the change of $\Delta(t)$. A scenario of rate-induced tracking requires that \eqref{eq:concave-variable-rate} is in \textsc{Case} \hyperlink{Cc}{{\rm C}} for small values of $d$ and in \textsc{Case} \hyperlink{Ac}{{\rm A}} for big values of $d$ (see also Section \ref{subsec:5Crateinduced}).

Next, we show how the features of $\Gamma,\Delta$ and $I$ can be chosen to guarantee a rate-induced tracking scenario. Assume that \eqref{eq:concave_logistic} admits $0<c_-<c_0<c_+$ such that $\boldsymbol\lambda_*(c_-,0)$ and $\boldsymbol\lambda_*(c_+,0)$ are both strictly negative and $\boldsymbol\lambda_*(c_0,0)>0$. Note that the continuity of $\boldsymbol\lambda_*(\cdot,0)$ provides the existence of at least two tipping points for
\begin{equation}\label{eq:5Bstaticc}
    x'=-\big(x-\Gamma(c\,t)\big)^2+I(t)
\end{equation}
in $(c_-,c_+)$.
Since $\Delta$ is continuous there is an open interval $(t_1,t_2)\subset \R$ such that $\boldsymbol\lambda_*(\Delta(t),0)>0$ for all $t\in(t_1,t_2)$. We shall assume that $0\in(t_1,t_2)$ and up to changing $c_0$ that $\Delta(0)=c_0$.
The continuous variation of solutions and the coercivity condition \hyperlink{C2}{{\rm\textbf{C2}}} guarantee that there is $d_0>0$ such that \eqref{eq:concave-variable-rate} is in \textsc{Case} \hyperlink{Cc}{{\rm C}} for all $d\in(0,d_0)$ because \eqref{eq:concave-variable-rate} can be considered as a small perturbation of \eqref{eq:5Bstaticc}$_{c_0}$, in a suitable compact interval containing $0$. 
The first ingredient for rate-induced tracking is achieved. 
In order to guarantee the second one (tracking for high rates) we shall use a Carathéodory analysis of switching at infinite rate on the line of the techniques in Ref.~\onlinecite{lno2}. 

We say that a transition equation \eqref{eq:concave-variable-rate} switches (or transitions at infinite rate) at time $t_0\in\R$ between the asymptotic problems \eqref{eq:5Bstaticc}$_{c_-}$ and \eqref{eq:5Bstaticc}$_{c_+}$, if $\Delta^{t_0}(t)=c_-$ for $t< t_0$ and $\Delta^{t_0}(t)=c_+$ for $t\ge t_0$. 
Note that, in this case, the locally pullback attractive solution $a_{\Delta^{t_0}}$ and the locally pullback repelling solution $r_{\Delta^{t_0}}$ (analogous of $a_\Gamma$ and $r_\Gamma$ of Subsection \ref{subsec:4Aconcavehypotheses}) are easily identifiable:
$a_{\Delta^{t_0}}(t)=a_{c_-}(t)$ for $t\le t_0$ and 
$r_{\Delta^{t_0}}(t)=r_{c_+}(t)$ for $t\ge t_0$, where $a_{c_-}$ is the locally pullback attractive solution of \eqref{eq:5Bstaticc}$_{c_-}$ and $r_{c_+}$ is the locally pullback repelling solution of \eqref{eq:5Bstaticc}$_{c_+}$.
The analysis of tipping and tracking for the switching problem is therefore extremely simple: \textsc{Case} \hyperlink{Ac}{{\rm A}} occurs if and only if $a_{c_-}(t_0)>r_{c_+}(t_0)$;
\textsc{Case} \hyperlink{Bc}{{\rm B}} occurs if and only if $a_{c_-}(t_0)=r_{c_+}(t_0)$;
\textsc{Case} \hyperlink{Cc}{{\rm C}} occurs if and only if $a_{c_-}(t_0)<r_{c_+}(t_0)$. Most importantly, reasoning as for Theorem 4.5 in Ref.~\onlinecite{lno2} one can show that if the switching problem is either in \textsc{Case} \hyperlink{Ac}{{\rm A}} or \textsc{Case} \hyperlink{Cc}{{\rm C}} then there is a minimum $\overline d>0$ such that \eqref{eq:concave-variable-rate} showcases the same dynamics for every $d> \overline d$. 
Therefore, the condition $a_{c_-}(0)>r_{c_+}(0)$ will provide the sought-for rate-induced tracking scenario.

\par\smallskip
\textbf{Time-dependent phase.} We now focus on a version of \eqref{eq:concave_logistic} with time-dependent phase, i.e.,
\begin{equation}\label{eq:concave_variable_phase}
x'=-\big[x-\Gamma\big(c(t-\Delta(dt))\big)\big]^2+I(t). 
\end{equation}
with $\Gamma:\R\to\R$ and $\Delta: \R\to[0,\infty)$ continuous functions respectively limiting to $\gamma_-$ and $s_-$ as $t\to-\infty$ and to $\gamma_+$ and $s_+$ as $t\to\infty$. The rate of the transition is now fixed at $c$, while the function $\Delta(t)$ represents a time-dependent phase and the new parameter $d>0$ is the "rate" for the change of $\Delta(t)$. Although now the function $\Delta$ has a different impact on the transition, much of the analysis shown in the previous paragraph can be carried out alike.
Particularly, values $0<s_-<s_0<s_+$ such that $\boldsymbol\lambda_*(c,s_-),\boldsymbol\lambda_*(c,s_+) <0$ and $\boldsymbol\lambda_*(c,s_0)>0$ guarantee that a value $d_c>0$ exists such that \eqref{eq:concave_variable_phase} is in \textsc{Case} \hyperlink{Cc}{{\rm C}} for all $d\in(0,d_c)$. 
Analogously, a suitably designed switching problem between $s_-$ and $s_+$ at time $t_0=0$ and the following Carathéodory analysis, provide a sufficient criterion $a_{s_-}(0)>r_{s_+}(0)$ to guarantee the occurrence of \textsc{Case} \hyperlink{Ac}{{\rm A}} for sufficiently high values of $d$ and thus rate-induced tracking.\par\smallskip

\textbf{Numerical evidence.} As follows, we aim to showcase the two previously discussed mechanisms generating rate-induced tracking for the concave logistic model \eqref{eq:concave_logistic} through numerical simulation. The results are shown in Figures \ref{fig:variable_rates} and \ref{fig:variable_phases}. In both cases the following choices of $\Gamma$, and $I$ have been used: 
\[
\begin{split}
    \Gamma(t)&=(2/\pi)\arctan(t)\\
    I(t)&=-\sin(t/2)-\sin(\sqrt5t)+0.895.
\end{split}
\]
In Figure \ref{fig:variable_rates} the case of time-dependent rate is shown. The variable transition $\Delta_{c_-,c_+}(t)=c_-/(1+e^t)+c_+/(1+e^{-t})$ between $c_-=0.25$ and $c_+=0.74$ assumes the value $c_0=(c_-+c_+)/2=0.495$ for $t=0$. As shown in the central panel of Figure \ref{fig:variable_rates} the approximation of the bifurcation map $\boldsymbol\lambda_*(c,0)$ for $c\in[0,1]$ is negative at $c_-$ and $c_+$ and positive at $c_0$. Further details on the approximation of $\boldsymbol\lambda_*(c,0)$ can be found in Ref.~\onlinecite{lno1}. The right-hand side panel of Figure \ref{fig:variable_rates} contains the approximation of the attractor $a_{c_-}$ and the repeller $r_{c_+}$ whose relative position at time $t=0$ determines the dynamic scenario of the switching problem. As explained in the previous paragraph, this evidence is sufficient to guarantee the phenomenon of rate-induced tracking for \eqref{eq:concave_variable_phase} as showcased in the four panels on the left-hand side of Figure \ref{fig:variable_rates} where also the role of FTLE as EWS is shown. \par \smallskip

The case of time-dependent phase is shown in Figure \ref{fig:variable_phases}, where the variable transition $\Delta_{s_-,s_+}(t)=s_-/(1+e^t)+s_+/(1+e^{-t})$ between $s_-=-5$ and $s_+=10$ takes the value $s_0=(s_-+s_+)/2=2.5$ at $t=0$. The analogous commentaries hold.

\subsection{Safe operating margins with time-dependent rate}\label{subsec:6Cnoreturnrateinducetracking}
In this subsection, we provide an example of mechanism acting as an immediate reaction to an EWS which can be used to avoid tipping in a rate-induced tracking scenario with time-dependent rate.
\begin{figure}
\centering
\includegraphics[width=0.45\textwidth]{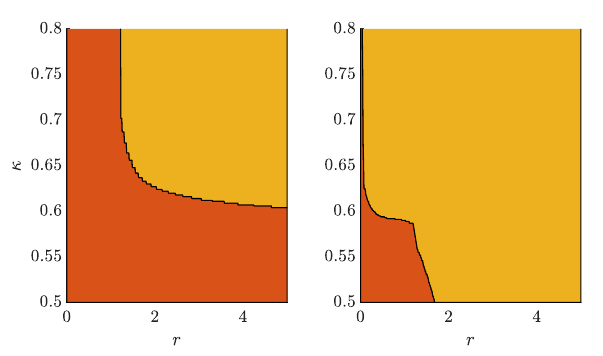}\par
\caption{Approximation of the regions of tipping and tracking for different reactions \eqref{eq:6Creaction} to early-warning signals triggered by an overgrowth of finite-time Lyapunov exponents. The functions $f$ and $\Gamma$ are chosen as in Figure~\ref{fig:6Cillustration}, and $b=1$. On the left side, $\Delta(t)=20-(\arctan(t/10)/\pi+1/2)$. On the right side, $\Delta(t)=20-0.75\cdot(\arctan(t/10)/\pi+1/2)$. In the $x$ axis, the strength $r$ of the reaction explained in \eqref{eq:6Creaction}. In the $y$ axis, the detection threshold $\kappa$: the bigger it is the earlier we detect early-warning signals. In yellow (above right), the region of pairs $(r,\kappa)$ for which the equation which takes the form \eqref{eq:6Cnonconstant} for $t\leq t_1$ and \eqref{eq:6Creaction} for $t\geq t_1$, where $t_1\in\mathbb{R}$ is the smallest value of time for which \eqref{eq:6Cwarning}, is in \textsc{Case} \protect\hyperlink{Ad}{{\rm A}}, in orange (below left), the region of pairs $(r,\kappa)$ for which such equation is in \textsc{Case} \protect\hyperlink{C2d}{{\rm C2}}.}\label{fig:6Cdetectionregions}
\end{figure}
Figure~\ref{fig:6Cdetectionregions} provides a numerical example on which we take specific action to reverse the tipping trend of an equation of the form \eqref{eq:6Cnonconstant} once the early-warning signal of finite-time Lyapunov exponents is detected.
For the sake of simplicity, we take a strictly decreasing $\Delta$ with asymptotic limits which crosses the rate-induced tracking tipping point $c_0$.
We assume that $|\Delta'|$ is small to ensure that the system inherits the early-warning signals of the constant $c$ problems.
Notice that this choice of $\Delta$ does not satisfy the requirements of Subsection~\ref{subsubsec:VA2geometricingredients}.
A detection threshold $\kappa\in[0,1)$ is taken as in Subsection \ref{subsec:5Crateinduced}, and $t_1\in\mathbb{R}$ is taken to be the smallest time where the early-warning signal is triggered,
\begin{equation}\label{eq:6Cwarning}
\frac{1}{T}\int_{t_1-T}^{t_1}f_x(s,x_\Delta(s),\Gamma(\Delta(s)\,s))\, ds\geq\kappa\,L\,.
\end{equation}
In order to reverse the tipping trend, the rate is hence increased starting from time $t_1$, leveraging on the fact that we are in a rate-induced tracking problem.
The increase of rate is carried out  by integrating,
\begin{equation}\label{eq:6Creaction}
x'=f\big(t,x,\Gamma\big(\big(\Delta(t)+r\tanh(b(t-t_1))\big)\,t\big)\big)\,,
\end{equation}
instead of \eqref{eq:6Cnonconstant}, where $r>0$ represents the strength of the reaction, which is an increase on the rate of the transition.
As expected, we find that the earlier we detect the warning (the bigger is $\kappa$) and the stronger we react to the warning-signal, the more probable is to avoid tipping.
\begin{table*}
\caption{\label{tab:table3}Summary of most common symbols used throughout the paper.}
\begin{ruledtabular}
\begin{minipage}{0.48\textwidth}
\begin{tabular}{cl}
Symbol & Meaning \\ \hline
$h$ & parameter independent total growth rate\\ \hline
$t\mapsto x(t,s,x_0)$ & maximal solution of $x'=h(t,x)$\\ \hline
 $f$ & parameter ($\gamma$) dependent total growth rate\\ \hline
 $\Gamma,\Delta$ & parameter shift function\\ \hline
 $\gamma_\pm$ & asymptotic values of $\Gamma$ \\ \hline
 $c_*,c_\pm,s_\pm$ & asymptotic values of $\Delta$ \\ \hline
$x'=f(t,x,\Gamma(t))$ & transition equation of parameter shift $\Gamma$ \\ \hline
$t\mapsto x_\Gamma(t,s,x_0)$ & maximal solution of $x'=f(t,x,\Gamma(t))$\\ \hline
$a_\Gamma$, $a_\Delta$ & \makecell[l]{locally pullback attractive solution of a \\ concave transition equation for $\Gamma$ or $\Delta$}\\ \hline
 $r_\Gamma$, $r_\Delta$ & \makecell[l]{locally pullback repulsive solution of a \\ concave transition equation for $\Gamma$ or $\Delta$}\\ \hline
 $u_\Gamma$, $u_\Delta$ & \makecell[l]{upper locally pullback attractive solution of a\\d-concave transition equation for $\Gamma$ or $\Delta$}\\ \hline
 $m_\Gamma$, $m_\Delta$ & \makecell[l]{locally pullback repulsive solution of a\\d-concave transition equation for $\Gamma$ or $\Delta$}\\ \hline
 $l_\Gamma$, $l_\Delta$ & \makecell[l]{lower locally pullback attractive solution of a\\d-concave transition equation for $\Gamma$ or $\Delta$}\\ \hline
 $x'=f(t,x,\gamma_\pm)$ & past and future equations of $x'=f(t,x,\Gamma(t))$\\ 
\end{tabular}
\end{minipage}\hfill\begin{minipage}{0.48\textwidth}
\begin{tabular}{cl}
Symbol & Meaning \\ \hline
  $\tilde a_{\gamma_\pm}$, $\tilde r_{\gamma_\pm}$& \makecell[l]{hyperbolic solutions of past and future equations\\ of a concave transition equation}\\ \hline
$\tilde u_{\gamma_\pm}$, $\tilde m_{\gamma_\pm}$, $\tilde l_{\gamma_\pm}$ & \makecell[l]{hyperbolic solutions of past and future equations\\ of a d-concave transition equation}\\ \hline
 $c,d$ & rate, phase or size parameter\\ \hline
 $\Gamma^c$ & parameter ($c$) dependent parameter shift function\\ \hline
 $x'=f(t,x,\Gamma^c(t))$ & parameter dependent transition equation\\  \hline
 $t\mapsto x_c(t,s,x_0)$ & maximal solution of $x'=f(t,x,\Gamma^c(t))$ \\ \hline
 $a_c$, $r_c$ & \makecell[l]{locally pullback attractive/repulsive solutions\\of a $c$-parametric concave transition equation}\\ \hline
 $u_c$, $m_c$, $l_c$ & \makecell[l]{locally pullback attractive/repulsive solutions\\of a d-concave $c$-parametric transition equation}\\ \hline
  $L$ & \makecell[l]{Lyapunov exponent of $\tilde a_{\gamma_-}$ (concave case),\\$\tilde u_{\gamma_-}$ or $\tilde l_{\gamma_-}$ (d-concave case)}\\ \hline
  \makecell[c]{$\lambda_a$, $\lambda_r$\\$\lambda_u$, $\lambda_m$, $\lambda_l$} & \makecell[l]{backward finite-time Lyapunov exponent of\\$a_c$, $r_c$, $u_c$, $m_c$ and $l_c$ respectively}\\\hline
      $\kappa\, L$ & warning threshold, with $\kappa\in[0,1)$\\ \hline
  $\boldsymbol\lambda_*(c,s)$ & \makecell[l]{bifurcation map for the concave equation \eqref{eq:concave_logistic}\\ with $\widehat\Gamma(t)=\Gamma(c(t+s))$} \\
\end{tabular}
\end{minipage}
\end{ruledtabular}
\end{table*}

\section{Conclusions}
A function $\Gamma$ with finite asymptotic limits can be understood as a parameter shift from its past asymptotic limit $\gamma_-$ to its future one $\gamma_+$.
If $\Gamma^c(t)$ is given by $\Gamma(ct)$ (for $c>0$) or $\Gamma(c+t)$ (for $c\in\mathbb R$), then each equation of the parametric family \eqref{eq:4Cparametrictransitionequation} can be understood as a transition from the past equation \eqref{eq:4pastequationgamma-} to the future equation \eqref{eq:4futureequationgamma+}.
We study critical transitions arising as $c$ changes for concave or d-concave (in $x$) functions $f$, assuming that the limit equations have the maximum number of hyperbolic solutions: two in the concave case, three in the d-concave case.
A critical transition occurs when the number of hyperbolic solutions of \eqref{eq:4Cparametrictransitionequation}$_c$ is the maximum for a value of the parameter, but no longer the maximum as $c$ crosses a tipping point $c_0$.
While the number of hyperbolic solutions is maximum, they connect the hyperbolic solutions of the past with those of the future, and it is said that \eqref{eq:4Cparametrictransitionequation}$_c$ shows \textsc{Case} \hyperlink{Ac}{{\rm A}} or tracking.
Typically, the tipping point $c_0$ corresponds to a nonautonomous saddle-node bifurcation point of the parametric family (see Section~\ref{sec:4tippingtracking}), where, upon the variation of $c$, two hyperbolic solutions, one attractive and the other repulsive, collide into a single nonhyperbolic solution at the tipping point $c_0$. 

It is important to keep in mind that nonautonomous saddle-node bifurcations can also admit more complicated dynamics at the tipping point. For example, a bifurcation point for  \eqref{eq:4Dlambdaparamteric} (whose existence can be proved under coercivity and concavity or d-concavity properties on $h$) does not necessarily imply the collision of two hyperbolic solutions into a nonhyperbolic one: their limits as $c\to c_0$ may be two different solutions (although they are no longer uniformly separated). 
However, the existence of the limit equations \eqref{eq:4pastequationgamma-} and \eqref{eq:4futureequationgamma+} with maximal number of hyperbolic solutions guarantees that a nonautonomous saddle-node bifurcation for \eqref{eq:x'=f(t,x,Gammac(t))} always entails that two hyperbolic solutions limit uniformly at the same unique nonhyperbolic solution as $c\to c_0$.  
In the skewproduct language (see Appendix \ref{appendix1}), this means that our equation corresponds to a point in its hull for which the two semicontinuous equilibria obtained as limit of the hyperbolic continuous equilibria take the same value (which happens at least at a residual set of points in the quasiperiodic case).
As a result of this property, the finite-time Lyapunov exponents (FTLE) calculated along a locally pullback attractive solution of the transition equation are efficient early-warning signals (EWS) alerting to the proximity of a tipping point: as the parameter $c$ approaches the critical value, the sign of the FTLE changes from negative to positive over time.

Throughout this work, we extensively investigate the presence of rate-induced tracking in rate parameter problems, that is, with $\Gamma^c(t)=\Gamma(ct)$ and $c>0$. There are types of transitions for which a tipping point may
occur when the positive transition rate $c$ decreases. This is the case 
of some of the models that we analyze: there are two bounded open sets $I$ 
and $J$ such that \eqref{eq:4intermediate}$_\gamma$ has the maximum number of hyperbolic solutions 
when $\gamma\in I$ but less if $\gamma\in J$; and there exists $T > 0$ 
such that $\Gamma(t)\in I$ for $|t|>T$ and $\Gamma(t)\in J$ for $|t|<T$. So, if the 
transition is fast, then $\Gamma^c$ takes ``safe'' values (in $I$) for a long
period of time, which guarantees the maximum number of hyperbolic solutions also for the 
transition equation \eqref{eq:4Cparametrictransitionequation}$_c$. But if the rate is low, then $\Gamma^c(t)$ is ``dangerous'' (in $J$) for a long period of time, and hence a locally pullback attractive solution departing from a steady state of the past system may blow up in the concave case or approximate an undesired solution in the d-concave case.

Section \ref{sec:Vrateinducedtrackingphenomena} shows that this behaviour can appear in a wide range of physical applications.
The conclusion of our numerical analysis is that the FTLE are also good EWS for
this type of problems even for a non-necessarily small transition rate.
In Subsection \ref{subsec:6Cnoreturnrateinducetracking}, we analyze transitions with time-dependent rate in some scenarios of rate-induced tracking.
We incorporate control mechanisms into the system that produce an abrupt acceleration of the transition immediately after an EWS, in order to force the dynamics of the transition to remain in \textsc{Case} \hyperlink{Ac}{{\rm A}}. This is possible only if this control takes place sufficiently early.
In other words, if $I$ and $J$ have the previously stated meaning, from the moment that the parameter shift $\Gamma$ starts assuming values in the dangerous set $J$, we only have a short interval of time to react and steer its values to the safe set of values $I$, avoiding thus tipping; but after this short interval of time to take action, the critical evolution of the transition has no-return.
The applicability of this result depends on the time scales involved in the problem: for instance, if time is measured in millennia, a ``short'' interval of time may correspond to enough years to incorporate those control mechanisms. Some of the issues mentioned in this final section will be the subject of further research in future publications.
\begin{acknowledgments}
The authors thank Carmen N\'{u}\~{n}ez for reading the manuscript and for her valuable suggestions that have improved the presentation of the work.

All the authors were supported by Ministerio de Ciencia, Innovaci\'{o}n y Universidades (Spain) under project PID2021-125446NB-I00 and by Universidad de Valladolid under project PIP-TCESC-2020. J. Due\~{n}as was supported by Ministerio de Universidades (Spain) under programme FPU20/01627. I.P.~Longo was partly supported by UKRI under the grant agreement EP/X027651/1, by the Horizon Europe - Societal Challenges project TiPES under grant agreement No. 820970, and by TUM International Graduate School of Science and Engineering (IGSSE) 
\end{acknowledgments}
\section*{Author declarations}
\vspace{-2ex}\noindent\textbf{Conflict of Interest}\par
The authors have no conflicts to disclose.\par\vspace{2ex}
\noindent\textbf{Author contributions}\par
\noindent\textbf{J.~Due\~{n}as:} Conceptualization, writing -- review \& editing.
\textbf{I.P.~Longo:} Conceptualization, writing -- review \& editing.
\textbf{R.~Obaya:} Conceptualization, writing -- review \& editing.
\section*{Data availability}
Data sharing is not applicable to this article as no new data were created or analyzed in this study.
\appendix
\section{Skewproduct formulation and Lyapunov exponents}\label{appendix1}
In this appendix, we present the skewproduct formulation, which is needed in the proofs of some of the results of Section~\ref{sec:4tippingtracking} as well as to define the Lyapunov spectrum of a bounded solution, which has been used in Section~\ref{subsec:2CLyapunovexponents}.
\subsection{Skewproduct formulation}\label{subsec:2Bskewproductpreliminaries}
The skewproduct formulation is a standard technique for the study of nonautonomous equations which allows us to use tools of topological dynamics.
Given $h\colon\mathbb{R}\times\mathbb{R}\rightarrow\mathbb{R}$, we represent by $h{\cdot}t(s,x)=h(t+s,x)$ the time translation of $h$ of magnitude $t\in\mathbb{R}$. Let $\mathcal F_h=\{h{\cdot}t\colon\, t\in\mathbb{R}\}$ be the set of all time translations of $h$, and let $\Omega_h$ be the closure of $\mathcal{F}_h$ on the space $C(\mathbb{R}\times\mathbb{R},\mathbb{R})$ endowed with the compact-open topology, that is, $\Omega_h$ contains not only the time translations of $h$ but also all the continuous functions which can be uniformly approximated on compact sets by time translations of $h$.
The set $\Omega_h$ is called the \emph{hull} of the function $h$.

It is known that, if $h\in C^{0,1}$ (resp. $h\in C^{0,2})$, then $\Omega_h$ is a compact metric subspace of $C^{0,1}$ (resp. $C^{0,2}$), that the time translation flow on the hull $\sigma\colon\mathbb{R}\times\Omega_h\rightarrow\Omega_h$, $\sigma(t,\omega)=\omega{\cdot}t$ is continuous, and that there exist continuous maps $\textswab{h},\textswab{h}_x\colon\Omega_h\times\mathbb{R}\rightarrow\mathbb{R}$ (and also $\textswab{h}_{xx}$ if $h\in C^{0,2}$) satisfying $\textswab h(\omega,x)=\omega(0,x)$ and $\textswab h_x(\omega,x)=\omega_x(0,x)$ for all $\omega\in\Omega_h$ and $x\in\mathbb{R}$.\cite{selltopdyn,shenyi}
If $h$ is quasiperiodic in time, then $\Omega$ is a minimal set (i.e., the closure of any of its $\sigma$-orbits) and supports a unique invariant (and hence ergodic) measure.

The skewproduct formulation of the problem consists on using the structure of the the hull $\Omega_h$ of $h$ given by $\sigma$ to deal jointly with all the equations of the family
\begin{equation}\label{eq:2Bfamilyskewproduct}
    x'=\omega(t,x)\,,\quad\text{for } \omega\in\Omega_h\,.
\end{equation}
For convenience, the family \eqref{eq:2Bfamilyskewproduct} is frequently rewritten by means of the previously introduced map $\textswab{h}$:
\begin{equation}\label{eq:2Bskewproductfamilyomega}
    x'=\textswab{h}(\omega{\cdot}t,x)\,,\quad\text{for }\omega\in\Omega_h\,.
\end{equation}
We refer to a equation of family \eqref{eq:2Bskewproductfamilyomega} for a particular fixed function $\omega\in\Omega_h$ as \eqref{eq:2Bskewproductfamilyomega}$_\omega$.
The family \eqref{eq:2Bskewproductfamilyomega} induces a (possibly local) \emph{skewproduct flow} on $\Omega_h\times\mathbb{R}$ given by $\tau\colon\mathcal{U}\subseteq\mathbb{R}\times\Omega_h\times\mathbb{R}\rightarrow\Omega_h\times\mathbb{R}$, $(t,\omega,x_0)=(\omega{\cdot}t,x(t,\omega,x_0))$, where $t\mapsto x(t,\omega,x)$ denotes the maximal solution of \eqref{eq:2Bskewproductfamilyomega}$_\omega$ with $x(0,\omega,x_0)=x_0$, that is, the maximal solution of \eqref{eq:2Bfamilyskewproduct}$_\omega$ for initial data $x_0$ at time $t=0$.

\subsection{Lyapunov exponents}\label{subsec:2CLyapunovexponents}
Hereafter, we will manage notions of orbits, invariant sets, ergodic measures, $\boldsymbol\alpha$-limit sets and $\boldsymbol\omega$-limit sets of a given flow.
The required definitions can be found in Ref.~\onlinecite{dno1}.
In what follows, $\mathfrak{M}_\mathrm{inv}(\Omega_h,\sigma)$ and $\mathfrak{M}_\mathrm{erg}(\Omega_h,\sigma)$ respectively stand for the nonempty sets of all normalized invariant and ergodic measures on $(\Omega_h,\sigma)$.

Let $\mathcal{K}\subset\Omega_h\times\mathbb{R}$ be a compact $\tau$-invariant set, that is, a compact set which is composed by orbits of the skewproduct flow $\tau$, and assume that $\mathcal{K}$ projects onto $\Omega_h$.
The \emph{Lyapunov exponent} on $\mathcal{K}$ for \eqref{eq:2Bskewproductfamilyomega} with respect to the invariant measure $\nu\in\mathfrak{M}_\mathrm{inv}(\mathcal{K},\tau)$ is
\begin{equation*}
\lambda_h(\mathcal{K},\nu)=\int_\mathcal{K} \textswab{h}_x(\omega,x)\, d\nu\,.
\end{equation*}
We will omit the subindex $h$ in $\lambda_h$ if there is no risk of confusion.
Theorems~1.8.4 of Ref.~\onlinecite{arnol} and 4.1 of Ref.~\onlinecite{furstenberg1} show that the Lyapunov exponent on $\mathcal{K}$ for \eqref{eq:2Bskewproductfamilyomega} with respect to any ergodic measure of $\mathfrak{M}_\mathrm{erg}(\mathcal{K},\tau)$ which projects onto an ergodic measure $m\in\mathfrak{M}_\mathrm{erg}(\Omega_h,\sigma)$ is given by an integral of the form 
\begin{equation}\label{eq:2Clyapunovexponents}
\int_\Omega \textswab{h}_x(\omega,\beta(\omega))\, dm\,,
\end{equation}
where $\beta\colon\Omega_h\rightarrow\mathbb{R}$ is an $m$-measurable equilibrium, that is, an $m$-measurable map with $\tau$-invariant graph.
The converse also holds: every integral of the form \eqref{eq:2Clyapunovexponents} with the graph of $\beta$ contained on $\mathcal{K}$ is a Lyapunov exponent on $\mathcal{K}$ for \eqref{eq:2Bskewproductfamilyomega}.
The \emph{Lyapunov spectrum} is the set of all Lyapunov exponents on $\mathcal{K}$ for \eqref{eq:2Bskewproductfamilyomega} with respect to all invariant measures of $\mathfrak{M}_\mathrm{inv}(\mathcal{K},\tau)$.
It is known that the Lyapunov spectrum on $\mathcal{K}$ for \eqref{eq:2Bskewproductfamilyomega} is a closed real interval whose endpoints are Lyapunov exponents with respect to ergodic measures.\cite{sackersell3,jps}
We can associate a \emph{Lyapunov spectrum} to any bounded solution $\tilde x\colon\mathbb{R}\rightarrow\mathbb{R}$ of \eqref{eq:2Ax'=h} by considering the Lyapunov spectrum on the smallest compact $\tau$-invariant set containing the $\tau$-orbit of $(h,\tilde x(0))$, that is, the set $\mathrm{cl}_{\Omega_h\times\mathbb{R}}\{(h{\cdot}t,\tilde x(t))\colon\, t\in\mathbb{R}\}$.
In the quasiperiodic case, this closure reduces to the graph of a continuous map if $\tilde x$ is hyperbolic, and hence the Lyapunov spectrum of any hyperbolic solution reduces to a point, which is negative if $\tilde x$ is attractive and positive if $\tilde x$ is repulsive.
\section{About the bifurcation results}\label{appendix2}
In this appendix, we give an approach to the proof of Theorem \ref{th:4Dd-concavebifurcationtheorem} to clarify the details which are not contained in the proof of Theorem 5.10 of Ref.~\onlinecite{dno1}. The main difference between these two results resides in the assumptions on the base $\Omega_h$ of the skewproduct flow, the hull of $h$.
In Theorem 5.10 of Ref.~\onlinecite{dno1}, $\Omega_h$ is assumed to be minimal, that is, every orbit is dense.
However, in Theorem \ref{th:4Dd-concavebifurcationtheorem}, we now assume that the orbit of $h$ is dense, but there can be other nondense orbits.
A flow is said to be \emph{transitive} if there exists a dense orbit, which is the case of $(\Omega_h,\sigma)$.
The next propositions provide a deeper understanding on the transitivity of $\Omega_h$.

\begin{proposition}\label{prop:Aresidualinvariantdense} The set $\Omega_0\subseteq\Omega_h$ of the points of $\Omega_h$ with dense orbit is a $\sigma$-invariant residual set.
\end{proposition}
\begin{proof} Analogous to that of Proposition~I.11.4 of Ref.~\onlinecite{mane1}.
%Since $\Omega_h$ is compact, given $n\in\mathbb{N}$, there exists a finite cover of $\Omega_h$ of open balls $\{B^n_j\}_{j=1,\dots,m_n}$ of radius $1/n$.
%Let us check that $\bigcap_{n=1}^\infty \mathcal{R}_n$ is a residual set of points with dense orbit, where $\mathcal{R}_n=\bigcap_{j=1}^{m_n}\bigcup_{t\in\mathbb{R}}\sigma_t(B_j^n)$.
%Taking into account the existence of a dense orbit in $\Omega_h$, it is not difficult to check that, for any $n\in\mathbb{N}$ and $j\in\{1,2,\dots,m_n\}$, the set $\bigcup_{t\in\mathbb{R}}\sigma_t(B_j^n)$ is residual.
%Therefore, $\bigcap_{n=1}^\infty \mathcal{R}_n$ is a residual set.
%Now, let us check that the orbit of any $\omega\in\bigcap_{n=1}^\infty\mathcal{R}_n$ is dense on $\Omega_h$.
%Given $\omega_0\in\Omega_h$ and $n\in\mathbb{N}$, there exist $j_n\in\{1,2,\dots,m_n\}$ such that $\omega_0\in B_{j_n}^n$, and $t_n\in\mathbb{R}$ such that $\sigma_{t_n}(\omega)\in B_{j_n}^n$, so the orbit of $\omega$ is dense. Since $\bigcap_{n=1}^\infty \mathcal{R}_n\subseteq\Omega_0$, the set $\Omega_0$ is residual.
\end{proof}
%The following proposition ensures that, if the transition is not contained in the past or in the future, then the residual set $\Omega_0$ corresponds to the orbit of any of the points.
\begin{proposition}\label{prop:AstructureofOmegatransitive} Let $\omega_0\in\Omega_0$, and let $\Omega_-$ be the $\boldsymbol\alpha$-limit of $\omega_0$ and $\Omega_+$ its $\boldsymbol\omega$-limit.
Then, $\Omega_h=\Omega_-\cup\{\omega_0{\cdot}t\colon\, t\in\mathbb{R}\}\cup\Omega_+$. In addition, if $\Omega_h\neq\Omega_-$ and $\Omega_h\neq\Omega_+$, then $\Omega_0=\{\omega_0{\cdot}t\colon\, t\in\mathbb{R}\}$.
\end{proposition}
\begin{proof} The first assertion is easy to check. Let us now see that $\Omega_h\neq\Omega_+$ ensures that
$\Omega_0\cap\Omega_+$ is empty: otherwise the closed invariant set $\Omega_+$ contains a dense orbit, a contradiction. Similarly, $\Omega_0\cap\Omega_-$ is empty if $\Omega_h\neq\Omega_-$, and hence $\Omega_0=\{\omega_0{\cdot}t\colon t\in\mathbb R\}$.
\end{proof}
Now, we will give some details of the proof of Theorem \ref{th:4Dd-concavebifurcationtheorem}. To this end, we consider the bifurcation problem
\begin{equation}\label{eq:Aparametricproblem}
    x'=\textswab{h}(\omega{\cdot}t,x)+\gamma\,,
\end{equation}
where $\textswab h\colon\Omega_h\times\mathbb{R}\rightarrow\mathbb{R}$ is the skewproduct form of $h$ (see Subsection \ref{subsec:2Bskewproductpreliminaries}), $h\in C^{0,2}$ is strictly d-concave uniformly for $t\in\mathbb{R}$ and $\lim_{|x|\rightarrow\infty} h(t,x)/x=-\infty$ uniformly for $t\in\mathbb{R}$.
Let $\mathfrak{l}_\gamma,\mathfrak{u}_\gamma\colon\Omega_h\rightarrow\mathbb{R}$ be the bounds of the global attractor $\mathcal{A}_\gamma$ of the skewproduct flow $\tau_\gamma$ induced by \eqref{eq:Aparametricproblem}$_\gamma$, that is, $\mathcal{A}_\gamma=\{(\omega,x)\colon\, \mathfrak{l}_\gamma(\omega)\leq x\leq\mathfrak{u}_\gamma(\omega)\}$.
Note that $\gamma\mapsto \mathfrak{l}_\gamma(\omega),\mathfrak{u}_\gamma(\omega)$ are strictly increasing for all $\omega\in\Omega_h$ (see Theorem~5.5. of Ref.~\onlinecite{dno1}).
 A \emph{copy of the base} is the invariant graph of a continuous equilibrium $\mathfrak{u}\colon\Omega_h\rightarrow\mathbb{R}$.
%By adapting the arguments of Proposition~2.8 of Ref.~\onlinecite{lineardissipativescalar} (see the final argument in the proof of Theorem~4.2(ii) of Ref.~\onlinecite{dno1}), we prove that all the Lyapunov exponents of a $\tau$-minimal set projecting onto $\Omega$ are strictly negative (resp.~positive) if and only if it is an attractive (resp. repulsive) hyperbolic copy of the base.
Theorem \ref{th:4Dd-concavebifurcationtheorem} is a consequence of Theorem \ref{th:Askewdconcavebifurcationtheorem}.

\begin{theorem}[D-concave bifurcation]\label{th:Askewdconcavebifurcationtheorem} Let $h\in C^{0,2}$ be strictly d-concave uniformly for $t\in\mathbb{R}$ and assume that $\lim_{|x|\rightarrow\infty} h(t,x)/x=-\infty$ uniformly for $t\in\mathbb{R}$. Assume that there exist $\omega_0\in\Omega_0$ (the residual invariant set of Proposition~{\rm\ref{prop:Aresidualinvariantdense}}) and $\gamma_0\in\mathbb{R}$ such that \eqref{eq:Aparametricproblem}$_{\gamma_0}^{\omega_0}$ has three uniformly separated bounded solutions.
Then, there exists an interval $I=(\gamma_1,\gamma_2)$ with $\gamma_0\in I$ such that
\begin{enumerate}[label=\rm{(\roman*)}]
\item for every $\gamma\in I$, there exist three disjoint compact $\tau_\gamma$-invariant sets $\mathcal K_{\gamma}^{l}$, $\mathcal K_{\gamma}^{m}$, $\mathcal K_{\gamma}^{u}$ which are hyperbolic copies of the base, given by the graphs of $\mathfrak{l}_\gamma<\mathfrak{m}_\gamma<\mathfrak{u}_\gamma$.
In addition, $\mathcal K_{\gamma}^{m}$ is repulsive and $\mathcal K_{\gamma}^{l}$, $\mathcal K_{\gamma}^{u}$ are attractive, $\gamma\mapsto\mathfrak{m}_\gamma$ is strictly decreasing on $I$, and the asymptotic dynamics of each equation is described by Theorem~{\rm3.3} of Ref.~\onlinecite{dno3}.
\item Let $\mathfrak{m}_{\gamma_1}(\omega)=\lim_{\gamma\downarrow\gamma_1}\mathfrak{m}_\gamma(\omega)$ and $\mathfrak{m}_{\gamma_2}(\omega)=\lim_{\gamma\uparrow\gamma_2}\mathfrak{m}_\gamma(\omega)$. Then, $\inf_{t\in\mathbb{R}}(\mathfrak{u}_{\gamma_1}(\omega{\cdot}t)-\mathfrak{m}_{\gamma_1}(\omega{\cdot}t))=0$ and $\inf_{t\in\mathbb{R}}(\mathfrak{m}_{\gamma_2}(\omega{\cdot}t)-\mathfrak{l}_{\gamma_2}(\omega{\cdot}t))=0$ for all $\omega\in\Omega_0$.
\item For $\gamma\in[\gamma_2,\infty)$ (resp. $\gamma\in(-\infty,\gamma_1]$), the graph of $\mathfrak{u}_\gamma$ (resp. $\mathfrak{l}_\gamma$) is an attractive hyperbolic copy of the base, and, for $\gamma\in(\gamma_2,\infty)$ (resp. $\gamma\in(-\infty,\gamma_1)$) $\lim_{t\rightarrow\infty}(\mathfrak{u}_\gamma(\omega{\cdot}t)-\mathfrak{l}_\gamma(\omega{\cdot}t))=0$ for all $\omega\in\Omega_0$.
In addition, there exists $\gamma_*\geq0$ sufficiently large such that $\mathcal{A}_\gamma$ is an attractive hyperbolic copy of the base for all $\gamma\not\in[-\gamma_*,\gamma_*]$.
\end{enumerate}
\end{theorem}
In what follows we will sketch the proof.
Analogous arguments to that of Proposition 3.16 of Ref. \onlinecite{dno3} show that the hypothesis on strict d-concavity of Theorem 5.10 of Ref. \onlinecite{dno1} is satisfied.
If $(\omega_0,x_1)$, $(\omega_0,x_2)$ are initial datum of two uniformly separated bounded solutions of \eqref{eq:Aparametricproblem}$_{\gamma}^{\omega_0}$, then, taking limits through convergent subsequences, we get that \eqref{eq:Aparametricproblem}$_{\gamma}^{\omega}$ has two uniformly separated solutions for all $\omega\in\mathrm{cl}\{\omega_0{\cdot}t\colon t\in\mathbb{R}\}$.
Besides, it follows easily that, if \eqref{eq:Aparametricproblem}$_{\gamma}^{\omega_0}$ has three separated solutions for some $\omega_0\in\Omega_0$ and $\gamma\in\mathbb{R}$, then
\eqref{eq:Aparametricproblem}$_{\gamma}^{\omega}$ has three separated solutions for all $\omega\in\Omega$ and there exist three hyperbolic copies of $\Omega$.

Starting from the three hyperbolic copies of the base of \eqref{eq:Aparametricproblem}$_{\gamma_0}^{\omega_0}$, the proof of Theorem 5.10 of Ref.~\onlinecite{dno1} can be repeated to ensure the existence of a maximal finite interval $I=(\gamma_1,\gamma_2)$, with $\gamma_0\in I$, on which (i) takes place.
The assertions in (ii) are consequences of the maximality of $I$.

%The last assertion in (iii) follows from Theorem~5.5(iii) of Ref.~\onlinecite{dno1}, so it only remain to indicate how to prove the other ones.
Let us prove (iii). First, we fix $\gamma\in[\gamma_2,\infty)$. The set $\mathcal K^u_\gamma=\mathrm{cl}_{\Omega_h\times\mathbb{R}}\{(\omega,\mathfrak{u}_\gamma(\omega))\colon\omega\in\Omega_h\}$ is a $\tau_\gamma$-invariant compact set.
The continuity of $\mathfrak{u}_{\gamma_0}$ and the inequality $\mathfrak{u}_{\gamma_0}<\mathfrak{u}_\gamma$ yield $\mathfrak{m}_{\gamma_0} (\omega)<\mathfrak{u}_{\gamma_0}(\omega)\leq\inf\{x\colon(\omega,x)\in\mathcal K^u_\gamma\}$.
Proposition~5.9 of Ref.~\onlinecite{dno1} ensures that all the Lyapunov exponents of $\mathcal{K}^u_\gamma$ are strictly negative: they have strictly negative sum with those of the repulsive hyperbolic compact $\tau_{\gamma_0}$-invariant set $\mathcal{K}^m_{\gamma_0}$, which are strictly positive.
 Therefore, $\mathcal{K}^u_\gamma$ is a hyperbolic attractive copy of the base: see the final argument in the proof of Theorem~4.2(ii) of Ref.~\onlinecite{dno1}.
 An easy contradiction argument shows that $\mathfrak{u}_\gamma$ is continuous. Hence, its graph provides $\mathcal K^u_\gamma$, which shows the first assertion.

 Let us now fix $\gamma\in(\gamma_2,\infty)$ and $\omega\in\Omega_0$, and prove that
$\lim_{t\rightarrow\infty}(\mathfrak{u}_\gamma(\omega{\cdot}t)-\mathfrak{l}_\gamma(\omega{\cdot}t))=0$.
We begin by checking that there exists $s\in\mathbb{R}$ such that $\mathfrak{l}_\gamma(\omega{\cdot}s)>\mathfrak{m}_{\gamma_2}(\omega{\cdot}s)$.
Otherwise, $\mathfrak{l}_{\gamma_2}(\omega{\cdot}t)<\mathfrak{l}_\gamma(\omega{\cdot}t)\leq\mathfrak{m}_{\gamma_2}(\omega{\cdot}t)$ for all $t\in\mathbb{R}$.
Then, $m_{\gamma_2}(\omega{\cdot}t)-l_{\gamma_2}(\omega{\cdot}t)\geq l_\gamma(\omega{\cdot}t)-l_{\gamma_2}(\omega{\cdot}t)=x_\gamma(1,\omega{\cdot}(t-1),l_\gamma(\omega{\cdot}(t-1)))-x_{\gamma_2}(1,\omega{\cdot}(t-1),l_{\gamma_2}(\omega{\cdot}(t-1)))\geq x_\gamma(1,\omega{\cdot}(t-1),l_{\gamma_2}(\omega{\cdot}(t-1)))-x_{\gamma_2}(1,\omega{\cdot}(t-1),l_{\gamma_2}(\omega{\cdot}(t-1)))$, so there exists $\xi\in[\gamma_2,\gamma]$ such that
\begin{equation*}
\begin{split}
    &m_{\gamma_2}(\omega{\cdot}t)-l_{\gamma_2}(\omega{\cdot}t)\\
    &\quad\geq\left.(\partial x_\mu/\partial \mu)\right|_\xi (1,\omega{\cdot}(t-1),l_{\gamma_2}(\omega{\cdot}(t-1)))(\gamma-\gamma_2)\,.
\end{split}
\end{equation*}
for some $\xi\in[\gamma_2,\gamma]$.
Since (ii) ensures that $\inf_{t\in\mathbb{R}}(\mathfrak{m}_{\gamma_2}(\omega{\cdot}t)-\mathfrak{l}_{\gamma_2}(\omega{\cdot}t))=0$, it suffices to check that there exists $k>0$ such that $\left.(\partial x_\mu/\partial \mu)\right|_\xi (1,\omega{\cdot}(t-1),l_{\gamma_2}(\omega{\cdot}(t-1)))\geq k$ for all $\xi\in[\gamma,\gamma_2]$ and $t\in\mathbb{R}$ to obtain a contradiction.
We define $y_{\xi,\omega_*}(s)=\left.(\partial x_\mu/\partial\mu)\right|_\xi(s,\omega_*,l_{\gamma_2}(\omega_*))$ for $(s,\xi,\omega_*)\in[0,1]\times[\gamma_2,\gamma]\times\Omega$ and observe that $y'_{\xi,\omega_*}(s)=\mathfrak h_x(s,x_\xi(s,\omega_*,l_{\gamma_2}(\omega_*))) y_{\xi,\omega_*}(s)+1$ and that $y_{\xi,\omega_*}(0)=0$.
Since the absolute value of $\mathfrak h_x(s,x_\xi(s,\omega_*,l_{\gamma_2}(\omega_*)))$ is bounded for $(s,\xi,\omega_*)\in [0,1]\times[\gamma_2,\gamma]\times\Omega$, solving the linear differential equation, it is easy to check that there exists $k>0$ such that $y_{\xi,\omega_*}(1)>k$ for $\xi\in[\gamma_2,\gamma]$ and $\omega_*\in\Omega$.
This proves the assertion.
We call $\bar\omega=\omega{\cdot}s$, so $\mathfrak{l}_\gamma(\bar\omega)>\mathfrak{m}_{\gamma_2}(\bar\omega)$.

Hence, there exists $\hat\gamma\in I$ such that $\mathfrak{l}_\gamma(\bar\omega)>\mathfrak{m}_{\hat\gamma}(\bar\omega)$.
According to the description of the dynamics of \eqref{eq:Aparametricproblem}$_{\bar\omega}^{\hat\gamma}$ given by Theorem 3.3 of Ref.~\onlinecite{dno3},
$\lim_{t\rightarrow\infty}(x_{\hat\gamma}(t,\bar\omega,\mathfrak{l}_\gamma(\bar\omega))-
\mathfrak{u}_{\hat\gamma}(\bar\omega{\cdot}t))=0$, and a comparison argument shows that $x_{\hat\gamma}(t,\bar\omega,\mathfrak{l}_\gamma(\bar\omega))<x_{\gamma}(t,\bar\omega,\mathfrak{l}_\gamma(\bar\omega))=\mathfrak{l}_\gamma(\bar\omega{\cdot}t)$ for all $t>0$.
Therefore, there exists $\bar t>0$ and $\bar\gamma\in(\gamma_1,\gamma)$ such that $\mathfrak{l}_\gamma(\bar\omega{\cdot}\bar t)=\mathfrak{u}_{\bar\gamma}(\bar\omega{\cdot}\bar t)$.

Now we define $A=\{\xi\in[\bar\gamma,\gamma)\colon$ there exists $s_\xi\geq \bar t$ such that
$\mathfrak{l}_\gamma(\bar\omega{\cdot}s_\xi)=\mathfrak{u}_\xi(\bar\omega{\cdot}s_\xi)\}$ and observe that $\bar\gamma\in A$.
If $\xi\in A$, then $\mathfrak{l}_\gamma(\bar\omega{\cdot}t)=x_\gamma(t,\bar\omega{\cdot} s_\xi,\mathfrak{u}_{\xi}(\bar\omega{\cdot}s_\xi))>x_\xi(t,\bar\omega{\cdot}
s_\xi,\mathfrak{u}_{\xi}(\bar\omega{\cdot}s_\xi))=\mathfrak{u}_\xi(\bar\omega{\cdot}t)$ for all $t>s_\xi$.
Using $\mathfrak{l}_\gamma(\bar\omega{\cdot}\bar t)=\mathfrak{u}_{\bar\gamma}(\bar\omega{\cdot}\bar t)$, the continuity and monotonicity of $\gamma\mapsto\mathfrak{u}_\gamma$ on $(\gamma_1,\infty)$ and the previous fact, it is not hard to check that $A$ is convex and open in $[\bar\gamma,\gamma)$, so $A=[\bar\gamma,\sup A)\subseteq[\bar\gamma,\gamma)$.
Let us check that $\sup A=\gamma$, assuming for
contradiction that $\sup A < \gamma$.
Given $\delta_0>0$, we define a compact set $K\subset\mathbb{R}$ such that $\mathfrak{u}_\xi(\bar\omega)\pm\delta_0\subset K$ for all $\xi\in[\bar\gamma,\gamma]$ and take $\epsilon\in(0,\gamma-\sup A)$.
Then, there exists $\delta\in(0,\delta_0)$ such that $|\mathfrak h(\bar\omega{\cdot}t,x)-\mathfrak h(\bar\omega{\cdot}t,y)|\leq\epsilon$ for $x,y\in K$ and $t\in\mathbb{R}$ if $|x-y|\leq\delta$.
The robustness of the hyperbolicity (see Theorem~2.2 of Ref.~\onlinecite{dno3}) allows us to fix $\xi\in A$ such that $\inf_{t\in\mathbb{R}}(\mathfrak{u}_{\sup A}(\bar\omega{\cdot}t)-\mathfrak{u}_{\xi}(\bar\omega{\cdot}t))\leq\delta$.
Then, $\mathfrak{l}_\gamma(\bar\omega{\cdot}t)=x_\gamma(t,\bar\omega{\cdot}\bar t,\mathfrak{u}_{\bar\gamma}(\bar\omega{\cdot}\bar t))\in(\mathfrak{u}_{\xi}(\bar\omega{\cdot}t),\mathfrak{u}_{\sup A}(\bar\omega{\cdot}t))$
for all $t>s_{\xi}$: otherwise, if there exists $s_{\sup A}>s_\xi$ such that $\mathfrak{l}_\gamma(\bar\omega{\cdot}s_{\sup A})\geq\mathfrak{u}_{\sup A}(\bar\omega{\cdot}s_{\sup A}))$, then either $\sup A=\gamma$ or the definition of $\sup A$ is contradicted.
So, $|\mathfrak{u}_{\sup A}(\bar\omega{\cdot}t)-x_\gamma(t,\bar\omega{\cdot}\bar t,\mathfrak{u}_{\bar\gamma}(\bar\omega{\cdot}\bar t))|\leq\delta$.
Therefore, for $t>s_{\xi}$,
\begin{equation*}
\begin{split}
\mathfrak{l}_\gamma(\bar\omega{\cdot}t)&=x_\gamma(t,\bar\omega{\cdot}\bar t,\mathfrak{u}_{\bar\gamma}(\bar\omega{\cdot}\bar t))\\ &=x_\gamma(s_{\xi},\bar\omega{\cdot}\bar t,\mathfrak{u}_{\bar\gamma}(\bar\omega{\cdot}\bar t))\\&\qquad\quad+\int_{s_{\xi}}^t \big(\mathfrak h(\bar\omega{\cdot}s,x_\gamma(s,\bar\omega{\cdot}\bar t,\mathfrak{u}_{\bar\gamma}(\bar\omega{\cdot}\bar t)))+\gamma\big)\,ds\\
&\geq x_\gamma(s_{\xi},\bar\omega{\cdot}\bar t,\mathfrak{u}_{\bar\gamma}(\bar\omega{\cdot}\bar t))+(\gamma-\sup A-\epsilon)(t-s_\xi)\\&\qquad\quad+\int_{s_{\xi}}^t \big(\mathfrak h(\bar\omega{\cdot}s,\mathfrak{u}_{\sup A}(\bar\omega{\cdot}s))+\sup A\big)\,ds\\
&= \mathfrak{u}_{\sup A}(\bar\omega{\cdot}t)+x_\gamma(s_{\xi},\bar\omega{\cdot}\bar t,\mathfrak{u}_{\bar\gamma}(\bar\omega{\cdot}\bar t))-\mathfrak{u}_{\sup A}(\bar\omega{\cdot}s_{\xi})\\&\qquad\quad+(\gamma-\sup A-\epsilon)(t-s_{\xi})\,.
\end{split}
\end{equation*}
The right-hand term tends to $\infty$, which contradicts the boundedness of $\mathfrak{l}_\gamma$.
Therefore, $A=[\bar\gamma,\gamma)$.
Hence, given $\epsilon>0$, there exists $\xi\in A$ such that $\inf_{t\in\mathbb{R}}(\mathfrak{u}_\gamma(\omega{\cdot}t)-\mathfrak{u}_{\xi}(\omega{\cdot}t))\leq\epsilon$, so $|\mathfrak{u}_\gamma(\omega{\cdot}t)-\mathfrak{l}_\gamma(\omega{\cdot}t))|\leq \epsilon$ for all $t\geq s_\xi+s$, that is,
$\lim_{t\rightarrow\infty}(\mathfrak{u}_\gamma(\omega{\cdot}t)-\mathfrak{l}_\gamma(\omega{\cdot}t))=0$, as we wanted.

The proofs of the assertions for $\gamma\in(-\infty,\gamma_1]$ and $\gamma\in(-\infty,\gamma_1)$ are analogous.
The last assertion in (iii) follows from Theorem 5.5(iii) of Ref.~\onlinecite{dno1}.

\bibliography{aipsamp}% Produces the bibliography via BibTeX.

\end{document}